# NONPARAMETRIC ESTIMATION OF MIXING DENSITIES FOR DISCRETE DISTRIBUTIONS

By François Roueff[1] and Tobias Rydén[2]

*GET/Telecom Paris and Lund University*

By a mixture density is meant a density of the form $\pi_\mu(\cdot) = \int \pi_\theta(\cdot) \times \mu(d\theta)$, where $(\pi_\theta)_{\theta \in \Theta}$ is a family of probability densities and $\mu$ is a probability measure on $\Theta$. We consider the problem of identifying the unknown part of this model, the mixing distribution $\mu$, from a finite sample of independent observations from $\pi_\mu$. Assuming that the mixing distribution has a density function, we wish to estimate this density within appropriate function classes. A general approach is proposed and its scope of application is investigated in the case of discrete distributions. Mixtures of power series distributions are more specifically studied. Standard methods for density estimation, such as kernel estimators, are available in this context, and it has been shown that these methods are rate optimal or almost rate optimal in balls of various smoothness spaces. For instance, these results apply to mixtures of the Poisson distribution parameterized by its mean. Estimators based on orthogonal polynomial sequences have also been proposed and shown to achieve similar rates. The general approach of this paper extends and simplifies such results. For instance, it allows us to prove asymptotic minimax efficiency over certain smoothness classes of the above-mentioned polynomial estimator in the Poisson case. We also study discrete location mixtures, or discrete deconvolution, and mixtures of discrete uniform distributions.

**1. Introduction.** Let $(\mathcal{X}, \mathcal{F})$ be a measurable space and let $(\pi_\theta)_{\theta \in \Theta}$ be a parametric family of densities on $\mathcal{X}$ with respect to a common measure $\zeta$. The parameter $\theta$ is assumed to range over a set $\Theta \in \mathbb{B}(\mathbb{R}^d)$; here $d \geq 1$ and $\mathbb{B}(\cdot)$ denotes the Borel sets. For any probability measure $\mu$ on $(\Theta, \mathbb{B}(\Theta))$, the

Received March 2002; revised October 2004.
[1]Supported by the Swedish Foundation for Strategic Research.
[2]Supported by a grant from the Swedish Research Council.

*mixture density* $\pi_\mu$ is defined on $\mathcal{X}$ by

$$\pi_\mu(x) = \int_\Theta \pi_\theta(x)\mu(d\theta).$$

Here the family $(\pi_\theta)$ is called the *mixands* and $\mu$ is the *mixing distribution*. If $\mu$ has finite support, $\pi_\mu$ is called a finite mixture (density). Estimation of such mixtures from an i.i.d. sequence $(X_i)_{1\leq i \leq n}$ distributed according to $\pi_\mu$, with the aim of recovering the unknown support points, their weights and maybe also their number, has a long history and we refer to the monographs by McLachlan and Peel [18], Titterington, Smith and Makov [20] and Lindsay [15] for further reading. In the present paper we are interested in nonparametric estimation of the mixing distribution $\mu$. We will assume that each such distribution under consideration has a density, called the *mixing density* and denoted by $f$, with respect to a known reference (Radon) measure $\nu$ on $(\Theta, \mathbb{B}(\Theta))$.

The problem of estimating $f$ for mixtures of discrete distributions ($\mathcal{X}$ is discrete) has been investigated, for instance, by Zhang [23] and, for Poisson mixtures with $\nu$ being Lebesgue measure, by Hengartner [12]; see also references in these two articles. The estimators examined by these authors are of two sorts. Zhang [23] used a kernel density estimator and adapted it to the mixture setting to estimate $f$ pointwise. Hengartner [12] used a projection estimator based on orthogonal polynomials to obtain an estimator of $f$ as an element of $L^2[a,b]$, $0 \leq a < b < \infty$. Loh and Zhang [16] used the kernel estimator to derive estimators of $f$ in the two cases $f \in L^p[0,b]$ and $f \in L^p[0,\infty)$ with $1 \leq p \leq \infty$. The main results of these works are concerned with establishing rates of convergence of the estimators, depending on smoothness conditions assumed on the mixing density, and with establishing bounds on the achievable minimax rate for mixing densities within balls defined by similar smoothness conditions.

The results on both estimators were condensed and slightly generalized by Loh and Zhang [17], who also carried out a numerical study of their finite sample performance. A conclusion of their work is that, although both types of estimators achieve similar rates with similar smoothness conditions on the mixing density, projection estimators seem to behave much better for finite samples. As pointed out by Loh and Zhang [17], the rates being logarithmic, it is not surprising that identical rates do not imply similar performance for finite sample sizes.

Another important point of the works cited above is that, although the rates of the estimators are derived over a wide range of smoothness classes, minimax rate optimality is proved only for particular instances. For example, Hengartner [12] obtained the rate of the projection estimator over Sobolev classes with arbitrary index of smoothness, but proved this rate to be minimax optimal for integer indices only. Similar remarks apply to the results of Loh and Zhang [17], but for a family of ellipsoidal classes.

In this paper we develop a general framework for studying projection estimators, with the main focus on mixtures of discrete distributions. Let us denote by $\Pi$ the linear operator mapping a real function $h$ on $\mathcal{X}$ to a real function $\Pi h$ on $\Theta$, defined by

$$\Pi h(\theta) = \pi_\theta h = \int_\mathcal{X} h \pi_\theta \, d\zeta \qquad \text{for all } \theta \text{ in } \Theta, \tag{1}$$

whenever this integral is well defined. Here we use the classical functional analysis notation $\pi h := \int h \, d\pi$. Above we defined $\pi_\theta$ and $\pi_\mu$ as densities on $\mathcal{X}$ with dominating measure $\zeta$, but we will also use the same notation for the corresponding probability measures. Observe that, by Fubini's theorem, for all $h$ such that $\pi_\mu |h| < \infty$,

$$\pi_\mu h = \int h(x) \left( \int_\Theta \pi_\theta(x) \mu(d\theta) \right) \zeta(dx) = \mu \Pi h. \tag{2}$$

The mean $\pi_\mu h$ may be estimated by a sample mean obtained using i.i.d. observations from $\pi_\mu$; see also [1], where this problem is addressed for $h$ within a given class of functions. The basic idea of what we call the *projection estimator* is now to estimate $\pi_\mu h$ for a suitable finite collection of functions $h$ and then to use (2) to obtain an estimate of $\mu$. The precise definition is given in Definition 1.

Our objective, classical in a nonparametric approach, is to find the asymptotic behavior of the minimax risk

$$\inf_{\hat{\mu}_n \in \mathcal{S}_n} \sup_{\mu \in \mathcal{C}} \pi_\mu^{\otimes n} l(\mu, \hat{\mu}_n),$$

where $\mathcal{C}$, $l$ and $\mathcal{S}_n$, respectively, denote a class of distributions, a loss function and a set of estimators defined on $\mathcal{X}^n$ and taking values in a set compatible with the choice of $l$; $\pi_\mu^{\otimes n}$ is the distribution of $n$ i.i.d. observations from $\pi_\mu$. It turns out that there is a simple argument to lower-bound this quantity in a general mixture framework (Proposition 1).

However, for exploiting this lower bound and studying the projection estimator, we will, as in the papers cited above, consider the case when $\mu$ is defined by its density $f = d\mu/d\nu$ for a fixed $\nu$. In this setting we will likewise write $\pi_f$ for $\pi_\mu$. Furthermore, the density $f$ will be assumed to belong to the Hilbert space $\mathbb{H} = L^2(\nu)$ with scalar product $(f, g)_\mathbb{H} = \int fg \, d\nu$. Given an estimator $\hat{f} : \mathcal{X}^n \to \mathbb{H}$ of $f$, it is natural to consider a risk given by the mean squared error $\mathbb{E}_f \| \hat{f} - f \|_\mathbb{H}^2$; here $\mathbb{E}_f$ denotes integration with respect to $\pi_f^{\otimes n}$ and $\| \cdot \|_\mathbb{H}$ is the norm on $\mathbb{H}$. In nonparametric language this is a mean integrated squared error (MISE). We will notice that, in order to arrive at interesting results, it is sensible to define the class $\mathcal{C}$ above, which is now a class of densities in $\mathbb{H}$, in accordance with the mixands. In the case of power series mixtures, this class is closely related to polynomials. Such ideas were used already by Lindsay [14] in a parametric framework. Still in the context of power series mixtures, we will obtain results on minimax rate optimality of the projection estimator, using classical results on polynomial approximations on compact sets (Theorem 3).

Having said that, we note that, quite generally, including Poisson mixands, the mixing density $f$ may also be estimated using nonparametric maximum likelihood; Lindsay [15] is excellent reading on this approach. The optimization problem so obtained is an infinite-dimensional convex programming problem, and numerical routines for approximating the nonparametric MLE (NPMLE) can be constructed, at least in certain models. The problem with the NPMLE is rather on the theoretical side. van de Geer [21] proved a rate of convergence result in terms of Hellinger distance in a rather abstract setting, and it still remains to be determined what this result implies for the problems studied in the present paper.

The paper is organized as follows. In Section 2 we give a general lower bound, in an abstract framework, on the obtainable error over certain classes of mixing distributions. This result is then specialized to the Hilbert setting outlined above, that is, we consider the MISE obtainable over smoothness classes of densities. In Section 3 we define the projection estimator and give a bias-variance decomposition of its loss. Section 4 focuses on mixtures of discrete distributions, containing a main theorem that provides a lower bound on the minimax MISE achievable over smoothness classes related to the definition of the projection estimator. An upper bound is also given and we discuss how these two bounds apply in a common setting. In Section 5 we apply these results to power series mixtures and complete the results obtained by Hengartner [12] and Loh and Zhang [17]. Section 6 is devoted to translation mixtures, or discrete deconvolution, while Section 7 provides applications of our results to mixtures of discrete uniform distributions. Finally, in Section 8 we give some examples in which the general methodology of the present paper may be valuable, but which we have not explored further.

Before closing this section we give some additional notation that will be used in connection with the above-mentioned Hilbert space $\mathbb{H}$. We write $\mathbb{H}_+$ for the set of nonnegative functions in $\mathbb{H}$, that is, $\mathbb{H}_+ = \{f \in \mathbb{H} : f \geq 0\}$, and $\mathbb{H}_1$ for the set of functions in $\mathbb{H}_+$ that integrate to unity, that is, $\mathbb{H}_1 = \{f \in \mathbb{H}_+ : \nu f = 1\}$. In other words, $\mathbb{H}_1$ is the set of probability densities on $\Theta$ which are also in $\mathbb{H}$. For any subset $V$ of $\mathbb{H}$, we write $V^\perp$ for the orthogonal complement of $V$ in $\mathbb{H}$, $f \perp V$ if $f \in V^\perp$ and $\mathrm{Proj}_V$ for the orthogonal projection on $V$. For two subsets $W \subseteq V$, we shall write $V \ominus W$ for $V \cap W^\perp$. A subset $V$ is called symmetric if $V = -V$, that is, if $-f$ is in $V$ whenever $f$ is.

**2. A general lower bound.** In this section we first give a lower bound on the obtainable loss in a more general framework, before turning to the setting specified in Section 1. To that end, let $(\mathcal{X}, \mathcal{F})$ and $(\Theta, \mathcal{G})$ be measurable spaces and let the function $(\theta, A) \mapsto \pi_\theta(A)$ from $\Theta \times \mathcal{G}$ to $[0, 1]$ be a probability kernel. That is, $\pi.(A)$ is measurable for all measurable subsets $A \subseteq \mathcal{X}$ and $\pi_\theta(\cdot)$ is a probability measure on $(\mathcal{X}, \mathcal{F})$ for all $\theta$ in $\Theta$. For instance, $\pi$ may be a regular version of a conditional probability $\mathbb{P}(X = \cdot | \theta = \cdot)$ (we refer to [19], Section 3.4, for more details).

We let $\mathbb{M}(\mathcal{X}, \mathcal{F})$ and $\mathbb{M}(\Theta, \mathcal{G})$ denote the sets of all signed finite measures on $(\mathcal{X}, \mathcal{F})$ and $(\Theta, \mathcal{G})$, respectively. For any set $A$, we write $\mathbf{1}_A$ for the indicator function of $A$. As in (1), the linear operator $\Pi$ maps a real function $h$ on $\mathcal{X}$ to a real function on $\Theta$ defined by $\Pi h(\theta) = \pi_\theta h$. Considering $\Pi$ acting on bounded functions, its adjoint operator $\Pi^*$ operates from $\mathbb{M}(\Theta, \mathcal{G})$ to $\mathbb{M}(\mathcal{X}, \mathcal{F})$ and is given by

$$\Pi^* \mu(A) = \Pi^* \mu \mathbf{1}_A = \mu \Pi \mathbf{1}_A \qquad \text{for all } A \in \mathcal{F}.$$

When $\mu$ is a probability measure on $\Theta$, its image $\Pi^* \mu$ by $\Pi^*$ also is a probability measure; indeed, it is the mixture distribution obtained from the mixands $(\pi_\theta)$ and mixing distribution $\mu$.

For any real function $h$ on $\mathcal{X}$, we denote by $M_h$ the multiplication operator which maps a real function $f$ on $\mathcal{X}$ to the real function defined by $M_h f(x) = h(x) f(x)$ for all $x \in \mathcal{X}$. Considering this operator acting on bounded functions, its adjoint $M_h^*$ is an operator on $\mathbb{M}(\mathcal{X}, \mathcal{F})$ given by $M_h^* \mu = \mu M_h$. In other words, $M_h^* \mu$ is the measure with density $h$ with respect to $\mu$.

Finally, we consider a subspace $E$ of signed measures on $\Theta$, equipped with a semi-norm $\mathcal{N}$. We assume that $E$ is endowed with a $\sigma$-field which makes this semi-norm measurable. We write $\mathcal{S}_n$ for the set of all $E$-valued estimators based on $n$ observations, that is, the set of all measurable functions from $\mathcal{X}^n$ to $E$. Finally, $E_1$ is the set of all probability measures that belong to $E$.

PROPOSITION 1. *Let $h$ be a real nonnegative function on $\mathcal{X}$, bounded by 1. Let $\mathcal{C}$ be a symmetric set included in the kernel of $M_h^* \circ \Pi^*$ and let $\mu_0$ be a probability measure on $\Theta$. Then for any number $p \geq 1$,*

$$\inf_{\hat{\mu} \in \mathcal{S}_n} \sup_{\mu \in (\mu_0 + \mathcal{C}) \cap E_1} (\Pi^* \mu)^{\otimes n} \mathcal{N}^p(\hat{\mu} - \mu)$$
(3)
$$\geq \sup\{\mathcal{N}^p(\mu) : \mu \in \mathcal{C}, \mu_0 \pm \mu \in E_1\} (\Pi^* \mu_0 h)^n.$$

REMARK. An obvious and interesting problem raised by the proposition is that of optimizing the right-hand side of (3) with respect to $h$.

REMARK. One can allow the function $h$ to depend on an index $i$ as long as $\mathcal{C}$ is a subset of the kernel of each $M_{h_i}^* \circ \Pi^*$. Doing so, the second factor in the lower bound becomes $\prod_{1 \leq i \leq n} \Pi^* \mu_0 h_i$.

REMARK. The supremum in the lower bound is also that of $\mathcal{N}^p$ over $\mathcal{C} \cap (E_1 - \mu_0) \cap (\mu_0 - E_1)$. Since this set is symmetric, the supremum is the $p$th power of its half diameter.

PROOF OF PROPOSITION 1. Write $h^{\otimes n}$ for the function on $\mathcal{X}^n$ mapping $(x_1, \ldots, x_n)$ to $\prod_{i=1}^n h(x_i)$. Let $\mu$ be in the kernel of $M_h^* \circ \Pi^*$ such that $\mu_0 + \mu$ is a probability measure. Then, since $h^{\otimes n}$ is bounded by 1, for all nonnegative functions $g$ on $\mathcal{X}^n$,

$$(\Pi^*(\mu_0 + \mu))^{\otimes n} g \geq (\Pi^*(\mu_0 + \mu))^{\otimes n}(g h^{\otimes n})$$
$$= (M_h^* \Pi^*(\mu_0 + \mu))^{\otimes n} g$$
$$= (M_h^* \Pi^* \mu_0)^{\otimes n} g.$$

Now pick an estimator $\hat{\mu}$ in $\mathcal{S}_n$. Let $\mu$ be a signed measure in the kernel of $M_h^* \circ \Pi^*$ such that both $\mu_\pm := \mu_0 \pm \mu$ are probability measures. Applying the above inequality twice, with $g$ equal to $g_\pm := \mathcal{N}^p(\hat{\mu} - \mu_\pm)$, we obtain

$$(4) \qquad (\Pi^* \mu_+)^{\otimes n} g_+ + (\Pi^* \mu_-)^{\otimes n} g_- \geq (M_h^* \Pi^* \mu_0)^{\otimes n}(g_+ + g_-).$$

Furthermore, note that

$$g_+ + g_- \geq 2^{1-p} \mathcal{N}^p(2\mu) = 2\mathcal{N}^p(\mu),$$

so that the right-hand side of (4) is at least $2\mathcal{N}^p(\mu)(\Pi^* \mu_0 h)^n$. The supremum in the left-hand side of (3) is at least half the left-hand side of (4), hence, at least $\mathcal{N}^p(\mu)(\Pi^* \mu_0 h)^n$. This corresponds to bounding the supremum risk from below by a two-point Bayes risk with uniform prior. We conclude the proof by optimizing over $\mu$. □

We note that Proposition 1 holds for norms such as the $L^p$ norms or the total variation norm in a nondominated context. Our particular interest in this result, however, is when $\pi_\theta(A) = \int_A \pi_\theta \, d\zeta$ (what we call the dominated case) and when $E$ is the set of all finite signed measures with a density with respect to $\nu$ in $\mathbb{H} = L^2(\nu)$. As this is the main topic of the remainder of the paper, we now restate Proposition 1 in this context as a separate result. From now on, $\mathcal{S}_n$ will denote the set of estimators in $\mathbb{H}$ from $n$ observations, that is, the set of measurable functions from $\mathcal{X}^n$ to $\mathbb{H}$, where $\mathbb{H}$ is endowed with its Borel $\sigma$-field.

PROPOSITION 2. *Let $f_0$ be in $\mathbb{H}_1$ and let $h$ be a real nonnegative function on $\mathcal{X}$, bounded by 1. Let $\mathcal{C}^\star$ be a symmetric subset of $\mathbb{H}$ such that, for $\zeta$-a.e. $x \in \mathcal{X}$, the mapping $\theta \mapsto h(x)\pi_\theta(x)$ belongs to $\mathcal{C}^{\star\perp}$. Then*

$$\inf_{\hat{f} \in \mathcal{S}_n} \sup_{f \in (f_0 + \mathcal{C}^\star) \cap \mathbb{H}_1} \mathbb{E}_f \|\hat{f} - f\|_\mathbb{H}^2$$
(5)
$$\geq \sup\{\|f\|_\mathbb{H}^2 : f \in \mathcal{C}^\star, f_0 \pm f \in \mathbb{H}_1\}(\pi_{f_0} h)^n.$$

PROOF. Take $E = \{M_f^*\nu : f \in \mathbb{H} \cap L^1(\nu)\}$ and define the norm $\mathcal{N}(M_f^*\nu) = \|f\|_\mathbb{H}$ on this space. Note that, for all $f$ in $\mathbb{H}_1$, $\Pi^* M_f^* \nu = \pi_f$. Thus, for all $f \in \mathbb{H}_1$, if $p = 2$ and $\mu = M_f^*\nu$, the expectation in the left-hand side of (1) equals that in the left-hand side of (5).

Now put $\mu_0 = M_{f_0}^*\nu$ and let $\mathcal{C} = \{M_f^*\nu : f \in \mathcal{C}^\star\} \cap \mathbb{M}(\Theta, \mathcal{G})$. From the assumptions on $\mathcal{C}^\star$, it is clear that $\mathcal{C}$ is a symmetric set included in the kernel of $M_h^* \circ \Pi^*$. Hence, in order to apply Proposition 1, it only remains to verify that

(6) $$\{M_f^*\nu : f \in (f_0 + \mathcal{C}^\star) \cap \mathbb{H}_1\} = (\mu_0 + \mathcal{C}) \cap E_1,$$

where $E_1$ is the set of probability distributions in $E$, that is, $E_1 = \{M_f^*\nu : f \in \mathbb{H}_1\}$. By observing that $\mathcal{C} \subseteq \{M_f^*\nu : f \in \mathcal{C}^\star\}$, we get the inclusion "$\supseteq$" in (6). For showing the inverse inclusion, pick $g \in \mathcal{C}^\star$ such that $f := f_0 + g$ is in $\mathbb{H}_1$. Since $f_0$ is in $\mathbb{H}_1$ as well, $M_g^*\nu \in \mathbb{M}(\Theta, \mathcal{G})$. This proves (6). □

Re-examining the proof of Proposition 1 in this context shows that it uses arguments similar to those of the second part of the proof of Theorem 3.1 in [12], but does not use Lemma 1 in [23], where a lower bound on $\sup_f \mathbb{P}_f\{\|\hat{f} - f\|_\mathbb{H} \geq \lambda\}$ is derived, the supremum being over a given subset of $\mathbb{H}_1$.

**3. The projection estimator.** Assume that $(X_i)_{1 \leq i \leq n}$ are i.i.d. with density $\pi_f$ with $f$ in $\mathbb{H}_1$. We denote by $P_n$ the empirical distribution defined by

$$P_n h = \int h \, dP_n = \frac{1}{n} \sum_{k=1}^n h(X_i) \qquad \text{for all } h : \mathcal{X} \to \mathbb{R}.$$

Let $\mathcal{H}$ denote the linear space containing all real functions $h$ which satisfy $\pi_\theta |h| < \infty$ for all $\theta$. For introducing the projection estimator, it is convenient to consider $\Pi$ defined by (1) acting on $\mathcal{H}$. The definition of the projection estimator depends on a given nondecreasing sequence $(V_m)_{m \geq 1}$ of finite-dimensional linear subspaces of $\mathbb{H}$. We put $d_m := \dim V_m$ and define $V_0 := \{0\}$. We assume without loss of generality that $V_m$ is included in $\Pi(\mathcal{H})$; otherwise we let $V_m \cap \Pi(\mathcal{H})$ replace $V_m$. We furthermore assume that, for any $g$ in $\bigcup_m V_m$, there exists a unique $h$ in $\mathcal{H}$ such that $\Pi h = g$, and we write $h = \Pi^{-1} g$. In other words, we assume that $\Pi$ is one-to-one on $\Pi^{-1}(\bigcup_m V_m)$. This is ensured, for instance, if $\Pi$ is one-to-one on $\mathcal{H}$, which, as observed by Barbe ([1], Lemma 5.1) simply means that the mixands are complete in the sense that, if $\Pi h(\theta) = 0$ for all $\theta$, then $h = 0$ (our $\Pi$ and $\mathcal{H}$ correspond to Barbe's $\mathcal{P}$ and $\mathcal{F}$, resp.). Moreover, he showed that for location and scale mixtures identifiability of the mixands in the sense $\pi_\mu = \pi_{\mu'}$ if and only if $\mu = \mu'$ implies that $\Pi$ is one-to-one ([1], Lemmas 5.2 and 5.3).

DEFINITION 1. Let $\widehat{f}_{m,n}$ be defined as the unique element in $V_m$ satisfying

$$(7) \qquad (\widehat{f}_{m,n}, g)_{\mathbb{H}} = P_n \Pi^{-1} g \qquad \text{for all } g \text{ in } V_m.$$

This estimator is called the *projection estimator* of $f$ of order $m$ [from $n$ observations and with respect to $(V_m)$].

From the assumptions above, the function which maps $g$ in $V_m$ to $P_n \Pi^{-1} g$ is a linear functional and, thus, (7) completely defines $\widehat{f}_{m,n}$ by duality of the scalar product. The projection estimator relies on the following idea. First observe that in the Hilbert setting (2) reads

$$(8) \qquad \pi_f h = (f, \Pi h)_{\mathbb{H}}$$

and holds for all $h$ such that $\Pi h$ is in $\mathbb{H}$. Hence, for all $g$ in $V_m$, by the law of large numbers, $P_n \Pi^{-1} g$ tends to $\pi_f \Pi^{-1} g = (f, \Pi \Pi^{-1} g)_{\mathbb{H}} = (f, g)_{\mathbb{H}}$ as $n \to \infty$, so that $\widehat{f}_{m,n}$ is approximately $\text{Proj}_{V_m} f$ for large $n$. Making $m$ large as well, $\text{Proj}_{V_m} f$ is roughly $f$, provided the closure of $\bigcup_{m \geq 1} V_m$ contains $f$. An important part of the development is thus to find a suitable rate at which to increase $m$ with respect to $n$.

In practice, the projection estimator can be expressed using an orthonormal sequence $(\phi_k)_{k \geq 0}$ in $\mathbb{H}$ such that $(\phi_k)_{0 \leq k \leq d_m - 1}$ is a basis of $V_m$ for all $m \geq 1$. The expansion of the projection estimator in this basis then reads

$$(9) \qquad \widehat{f}_{m,n} = \sum_{k=0}^{d_m - 1} (P_n \Pi^{-1} \phi_k) \phi_k.$$

For any random element $g$ in $\mathbb{H}$ such that $\pi_f \|g\|_{\mathbb{H}}^2 < \infty$, we define its variance as $\text{var}_f(g) := \mathbb{E}_f \|g - \mathbb{E}_f g\|_{\mathbb{H}}^2$. Under the i.i.d. assumption, the MISE of the projection estimator admits the following bias-variance decomposition.

PROPOSITION 3. *For all $f$ in $\mathbb{H}_1$, the MISE of $\widehat{f}_{m,n}$ writes*

$$(10) \qquad \mathbb{E}_f \|\widehat{f}_{m,n} - f\|_{\mathbb{H}}^2 = \|f - \text{Proj}_{V_m} f\|_{\mathbb{H}}^2 + \frac{1}{n} \text{var}_f(\widehat{f}_{m,1}).$$

PROOF. Pythagoras' theorem gives
$$\|\widehat{f}_{m,n} - f\|_{\mathbb{H}}^2 = \|\widehat{f}_{m,n} - \text{Proj}_{V_m} f\|_{\mathbb{H}}^2 + \|f - \text{Proj}_{V_m} f\|_{\mathbb{H}}^2.$$

From (9) and (8), we have
$$\mathbb{E}_f \widehat{f}_{m,n} = \sum_{k=0}^{d_m-1} \mathbb{E}_f (P_n \Pi^{-1} \phi_k) \phi_k = \sum_{k=0}^{d_m-1} (f, \phi_k)_{\mathbb{H}} \phi_k = \text{Proj}_{V_m} f.$$

Inserting this equality into the next to last display and taking expectations yields $\mathbb{E}_f \|\widehat{f}_{m,n} - f\|_{\mathbb{H}}^2 = \|f - \text{Proj}_{V_m} f\|_{\mathbb{H}}^2 + \text{var}_f(\widehat{f}_{m,n})$. Using (9) and the orthonormality of $(\phi_k)$, we obtain

$$(11) \quad \text{var}_f(\widehat{f}_{m,n}) = \sum_{k=0}^{d_m-1} \text{var}_f(P_n \Pi^{-1} \phi_k) = \frac{1}{n} \sum_{k=0}^{d_m-1} \text{var}_f(P_1 \Pi^{-1} \phi_k).$$

The proof is complete. □

We finish this section by noting that in many cases the sequence $(V_m)$ is defined as $V_m = \text{Span}(\Pi h_0, \ldots, \Pi h_{d_m-1})$ for a sequence $(h_k)_{k \geq 0}$ in $\mathcal{H}$ such that $(\Pi h_k)_{k \geq 0}$ is a linearly independent sequence in $\mathbb{H}$. This constructive definition of $(V_m)$ automatically ensures that all the above assumptions are verified. Observe, however, that the projection estimator only depends on the sequence $(V_m)$, whence different choices of $(h_k)$ are possible. In particular, by the Gram–Schmidt procedure, we can construct an orthonormal sequence $(\phi_k)$ as

$$\phi_k = \sum_{\ell=0}^{k} \Phi_{k,\ell} \Pi h_\ell \quad \text{for all } k \geq 0,$$

for some real coefficients $(\Phi_{k,\ell})_{k,\ell \geq 0}$ for which we set $\Phi_{k,\ell} := 0$ for all $\ell > k \geq 0$. The sequence $(\phi_k)$ may then replace $(h_k)$ for defining the same sequence $(V_m)$, and in this context (9) becomes

$$(12) \quad \widehat{f}_{m,n} = \sum_{k=0}^{d_m-1} \sum_{\ell=0}^{k} \Phi_{k,\ell} (P_n h_\ell) \phi_k.$$

**4. Application to mixtures of discrete distributions.** The basic assumption of this section is

$$\mathcal{X} = \mathbb{Z}_+ \text{ and } \zeta \text{ is counting measure.}$$

The case of continuous $\mathcal{X}$ seems to require deep adaptations and is left for future work. In the present setting we write $\mathbf{1}_k$ for the indicator function $\mathbf{1}_k(x) = \mathbf{1}(x = k)$ and take

$$(13) \quad V_m := \text{Span}(\Pi \mathbf{1}_k, 0 \leq k < m).$$

Notice that $\Pi \mathbf{1}_k = \pi.(k)$. We are hence in the constructive framework of Section 3 with $\dim V_m = m$, provided that $(\pi.(k))_{k \geq 0}$ is a sequence of linearly independent functions in $\mathbb{H}$. In this section we thus make the following assumption.

(A1) $(\Pi \mathbf{1}_k)_{k \geq 0}$ is a sequence of linearly independent functions in $\mathbb{H} \cap L^1(\nu)$.

Obviously, since $(\mathbf{1}_k)$ is a linearly independent sequence in $\mathcal{H}$, so is $(\Pi \mathbf{1}_k)$, provided $\Pi$ is one-to-one. We recall that this holds whenever the mixands are complete (see Section 3). Assumption (A1) implies that the projection estimator $\widehat{f}_{m,n}$ is well defined and, as a linear combination of $\Pi \mathbf{1}_k$'s, belongs to $L^1(\nu)$ for all $m$ and $n$. Hence, it is a good candidate for estimating a probability density function with respect to $\nu$. We elaborate further on this assumption in Section 4.3.

The results of Sections 2 and 3 may be used for bounding the minimax MISE $\inf_{\hat{f} \in \mathcal{S}_n} \sup_{f \in \mathcal{C}} \mathbb{E}_f \|\hat{f} - f\|_{\mathbb{H}}^2$ for particular smoothness classes $\mathcal{C}$, which we now introduce.

For any positive decreasing sequence $u = (u_m)_{m \geq 0}$, any positive number $C$ and any nonnegative integer $r$, define

(14) $\quad \mathcal{C}(u, C, r) := \{ f \in \mathbb{H} : \|f - \text{Proj}_{V_m} f\|_{\mathbb{H}} \leq C u_m \text{ for all } m \geq r \}.$

Note that, for $r \geq 1$, the classes $\mathcal{C}(u, 0, r)$ do not reduce to $\{0\}$ but to $V_r$. Also note that one may assume $u_0 = 1$ without loss of generality, in which case, recalling the convention $V_0 = \{0\}$,

(15) $\begin{aligned}&\mathcal{C}(u, C, 0) \\ &= \{ f \in \mathbb{H} : \|f\|_{\mathbb{H}} \leq C, \|f - \text{Proj}_{V_m} f\|_{\mathbb{H}} \leq C u_m \text{ for all } m \geq 1 \}.\end{aligned}$

Usually we simply write $\mathcal{C}(u, C)$ for $\mathcal{C}(u, C, 0)$. This set can be interpreted as the ball of functions whose rate of approximation by projections on the spaces $V_m$ is controlled by $(u_m)$ within a radius $C$. Finally, observe that having $\lim u_m = 0$ amounts to saying that $\mathcal{C}(u, C, r)$ is a subset of the closure of $\bigcup_{m \geq 1} V_m$ in $\mathbb{H}$.

For any fixed $f_0$ in $\mathbb{H}_+$, we define the following semi-norm on $\mathbb{H}$:

(16) $\quad \|f\|_{\infty, f_0} := \nu\text{-}\operatorname*{ess\,sup}_{\theta \in \Theta} \frac{|f(\theta)|}{f_0(\theta)},$

with the convention $0/0 = 0$ and $s/0 = \infty$ for $s > 0$. This semi-norm is not necessarily finite. Also introduce, for any subspace $V$ of $\mathbb{H}$,

$$K_{\infty, f_0}(V) := \sup\{ \|f\|_{\infty, f_0} : f \in V, \|f\|_{\mathbb{H}} = 1 \}.$$

Finally, we define for any positive numbers $K$ and $C$, any sequence $u = (u_m)$ as above and any nonnegative integer $r$,

$$\begin{aligned}
\mathcal{C}_{f_0}(K, u, C, r) &:= \{f \in \mathcal{C}(u, C, r) : \|f\|_{\infty, f_0} \leq K\} \\
&= \mathcal{C}(u, C, r) \cap \{f \in \mathbb{H} : \|f\|_{\infty, f_0} \leq K\};
\end{aligned} \tag{17}$$

again, just as for $\mathcal{C}(u, C)$, we write $\mathcal{C}_{f_0}(K, u, C)$ for $\mathcal{C}_{f_0}(K, u, C, 0)$.

4.1. *A lower bound on the MISE under* (A1). The following result is derived from Proposition 2 using the smoothness classes above.

THEOREM 1. *Let $f_0$ be in $\mathbb{H}_1$, $u = (u_m)_{m \geq 0}$ a positive decreasing sequence, $C$ a positive number, $r$ a nonnegative integer and $K$ a positive number such that $K \leq 1$. Then for any positive integer $n$, any estimator $\hat{f}_n$ in $\mathcal{S}_n$ and any integer $m \geq r$,*

$$\begin{aligned}
&\sup_{f \in (f_0 + \mathcal{C}_{f_0}(K, u, C, r)) \cap \mathbb{H}_1} \mathbb{E}_f \|\hat{f} - f\|_{\mathbb{H}}^2 \\
&\qquad \geq \left( \frac{K}{K_{\infty, f_0}(V_{m+2} \ominus V_m)} \wedge (C u_{m+1}) \right)^2 (\pi_{f_0}\{0, \ldots, m-1\})^n.
\end{aligned} \tag{18}$$

REMARK. For the lower bound (18) to be nontrivial, $K_{\infty, f_0}(V_{m+2} \ominus V_m)$ must be finite. Since $V_{m+2} \ominus V_m$ is finite-dimensional, this is true if $\|\cdot\|_{\infty, f_0}$ is a *finite* norm on $V_{m+2} \ominus V_m$.

REMARK. The lower bound (18) can be optimized over all $m \geq r$. In most cases $K_{\infty, f_0}(V_{m+2} \ominus V_m)$ behaves like $K_{\infty, f_0}(V_m)$ and thus increases as $m$ gets large. Hence, the squared term in the lower bound decreases when $m$ gets large while, in contrast, $\pi_{f_0}\{0, \ldots, m-1\}$ increases to 1 as $m$ tends to infinity.

The proof of the theorem is prefaced by two lemmas.

LEMMA 1. *Let $f_0$ be in $\mathbb{H}_+$. Then for all $f$ in $\mathbb{H}$,*

$$\|f\|_{\infty, f_0} = \sup_{g \in \mathbb{H}} \frac{|(f, g)_{\mathbb{H}}|}{(f_0, |g|)_{\mathbb{H}}} \tag{19}$$

*with the convention $0/0 = 0$ and $s/0 = \infty$ for $s > 0$.*

PROOF. First assume that there is a Borel subset $A$ of $\Theta$ with $\nu(A) > 0$ and such that both $f_0 = 0$ and $|f| > 0$ on $A$. It then follows immediately that the left-hand side of (19) is infinite, and so is the right-hand side (take $g = f \mathbf{1}_A$).

Now assume that there is no such set $A$. Using the convention $0/0 = 0$, we then have $f = (f/f_0)f_0$ $\nu$-a.e. Letting $\mu_0$ be the measure having density $f_0$ with respect to $\nu$, we find that the left-hand side of (19), $\nu\text{-}\operatorname{ess\,sup}|f/f_0|$, equals $\mu_0\text{-}\operatorname{ess\,sup}|f/f_0|$. Furthermore, $\mu_0$ is a $\sigma$-finite measure. Indeed, since $\nu$ is $\sigma$-finite, the Cauchy–Schwarz inequality shows that $\mu_0(K) = (f_0, \mathbf{1}_K)_{\mathbb{H}} \leq \|f_0\|_{\mathbb{H}}(\nu\mathbf{1}_K)^{1/2} < \infty$ for any compact set $K$. Hence, the space $L^\infty(\mu_0)$ and the dual $L^1(\mu_0)^*$ are isometric (see [9], Theorem 4.14.6), implying that the left-hand side of (19) equals

$$\mu_0\text{-}\operatorname{ess\,sup}|f/f_0| = \sup_{g:\,\mu_0|g|=1} |\mu_0[(f/f_0)g]| = \sup_{g:\,\mu_0|g|<\infty} \frac{|\mu_0[(f/f_0)g]|}{\mu_0|g|},$$

again with the convention $0/0 = 0$. It now remains to show that this display is equal to the right-hand side of (19).

To do that, notice that, for any $g$ in $\mathbb{H}$, $(f,g)_{\mathbb{H}} = \nu(fg) = \mu_0[(f/f_0)g]$ and $(f_0, |g|)_{\mathbb{H}} = \mu_0|g|$. Thus, the right-hand side of (19) is the supremum of the same ratio as in right-hand side of the last display, but over $g$ in $\mathbb{H}$ rather than over $g$ in $L^1(\mu_0)$. However, these suprema are, in fact, identical, which concludes the proof. To see the equality, first observe that since $\mu_0|g| = (|g|, f_0)_{\mathbb{H}} \leq \|g\|_{\mathbb{H}}\|f_0\|_{\mathbb{H}}$ for any $g$ in $\mathbb{H}$ (Cauchy–Schwarz), $\mathbb{H}$ is included in $L^1(\mu_0)$. The inverse inclusion does not hold, but, by optimizing the sign of $g$ in the two suprema, we may replace $f$ by $|f|$ in the numerators and restrict the suprema to nonnegative $g$'s and then use the result that any nonnegative function $g$ in $L^1(\mu_0)$ can be approximated by an increasing sequence of functions in $\mathbb{H}$ [e.g., by $(g\mathbf{1}_{|g|\leq M})_{M>0}$]. $\square$

LEMMA 2. *Adopt the assumptions of Theorem* 1 *and denote by* $\mathcal{C}^\star$ *the set* $\mathcal{C}_{f_0}(K,u,C,r) \cap V_m^\perp$. *We then have the upper and lower bounds*

(20) $$\sup\{\|f\|_{\mathbb{H}} : f \in \mathcal{C}^\star, f_0 \pm f \in \mathbb{H}_1\} \leq Cu_m$$

*and*

(21) $$\sup\{\|f\|_{\mathbb{H}} : f \in \mathcal{C}^\star, f_0 \pm f \in \mathbb{H}_1\} \geq Cu_{m+1} \wedge \frac{K}{K_{\infty,f_0}(V_{m+2} \ominus V_m)}.$$

PROOF. We start with the upper bound (20). Pick $f$ in $\mathcal{C}^\star$. Since $f$ is then in $\mathcal{C}(u,C,r)$, $\|f - \operatorname{Proj}_{V_m} f\|_{\mathbb{H}} \leq Cu_m$ for $m \geq r$. However, because $f$ is also in $V_m^\perp$, $\operatorname{Proj}_{V_m} f = 0$, and, thus $\|f\|_{\mathbb{H}} \leq Cu_m$.

We now turn to the lower bound. Let $(\phi_k)_{k\geq 0}$ be an orthonormal sequence in $\mathbb{H}$ such that $V_m = \operatorname{Span}(\phi_0, \ldots, \phi_{m-1})$ for all $m \geq 1$ (see Section 3). Using the fact that $\sum_{k\geq 0} \Pi \mathbf{1}_k = \Pi 1 = 1$, monotone convergence provides

$$\sum_{k\geq 0} (\Pi \mathbf{1}_\ell, \Pi \mathbf{1}_k)_{\mathbb{H}} = (\Pi \mathbf{1}_\ell, 1)_{\mathbb{H}} = \nu \Pi \mathbf{1}_\ell \qquad \text{for all } \ell \geq 0.$$

The right-hand side of this equation is finite by (A1). Since $\phi_l$ is a linear combination of $(\Pi \mathbf{1}_s)_{0 \leq s \leq l}$, we obtain

(22) $$\sum_{k \geq 0} |(\phi_\ell, \Pi \mathbf{1}_k)_\mathbb{H}| < \infty \qquad \text{for all } \ell \geq 0.$$

We shall now prove (21) by constructing a function $f$ in $\mathcal{C}^\star$ satisfying $f_0 \pm f \in \mathbb{H}_1$ and whose norm equals the right-hand side of (21). To that end, note that, by (22), we can find two numbers $\alpha$ and $\beta$ such that

(23) $$\alpha \sum_{k \geq 0} (\phi_m, \Pi \mathbf{1}_k)_\mathbb{H} + \beta \sum_{k \geq 0} (\phi_{m+1}, \Pi \mathbf{1}_k)_\mathbb{H} = 0$$

and, putting $f := \alpha \phi_m + \beta \phi_{m+1}$,

(24) $$\|f\|_\mathbb{H} = (\alpha^2 + \beta^2)^{1/2} = Cu_{m+1} \wedge \frac{K}{K_{\infty, f_0}(V_{m+2} \ominus V_m)}.$$

To finish the proof, we need to show that $f \in \mathcal{C}^\star$ and $f_0 \pm f \in \mathbb{H}_1$. To start with we note that $f$ lies in $V_{m+2}$ and that $f \perp V_m$. Therefore, $\|f - \operatorname{Proj}_{V_p} f\|_\mathbb{H} = 0$ for all $p \geq m+2$. Moreover, $\|f - \operatorname{Proj}_{V_{m+1}} f\|_\mathbb{H} = |\beta| \leq Cu_{m+1}$ and, since $(u_n)$ is decreasing, $\|f - \operatorname{Proj}_{V_p} f\|_\mathbb{H} = (\alpha^2 + \beta^2)^{1/2} \leq Cu_p$ for all $p = r, \ldots, m$. All this implies that $f$ lies in $\mathcal{C}(u, C, r)$. Using (24), we also see that $\|f\|_{\infty, f_0} \leq \|f\|_\mathbb{H} K_{\infty, f_0}(V_{m+2} \ominus V_m) \leq K$, so that $f$ belongs to $\mathcal{C}(K, u, C, r)$. Thus, $f \in \mathcal{C}^\star$.

Finally, as a finite linear combination of $L^1(\nu)$ functions, $f$ is in $L^1(\nu)$. Hence, dominated convergence and (23) yield $\nu f = \sum_{k \geq 0} (f, \Pi \mathbf{1}_k)_\mathbb{H} = 0$. By Lemma 1, we also find that, for all $g$ in $\mathbb{H}_+$,

$$(f_0 + f, g)_\mathbb{H} \geq (f_0, g)_\mathbb{H} (1 - \|f\|_{\infty, f_0}) \geq 0,$$

where we have used $K \leq 1$. Taking $g = (f_0 + f)_- := -(f_0 + f) \vee 0$ yields $-\|(f_0 + f)_-\|_\mathbb{H} \geq 0$, whence $(f_0 + f)_- = 0$ and $f_0 + f \in \mathbb{H}_+$. Together with $\nu f = 0$, this shows that $f_0 + f \in \mathbb{H}_1$. The same arguments hold true for $f_0 - f$ and the proof is complete. $\square$

PROOF OF THEOREM 1. Take $\mathcal{C}^\star$ as in Lemma 2 and define $h : \mathcal{X} \to \{0, 1\}$ by $h(x) := \mathbf{1}(0 \leq x < m)$. Then any mapping $\theta \mapsto h(x) \pi_\theta(x)$ is either identically zero (if $x \geq m$) or equal to $\pi_\theta(x) = \Pi \mathbf{1}_x$ (when $x < m$). Since such a $\Pi \mathbf{1}_x$ trivially lies in $V_m$, it is orthogonal to $\mathcal{C}^\star$. Thus, the conditions of Proposition 2 are met. Proposition 2, Lemma 2 and the trivial observation $\mathcal{C}^\star \subseteq \mathcal{C}_{f_0}(K, u, C, r)$ now prove the theorem. $\square$

4.2. *An upper bound on the MISE under* (A1). We shall now derive an upper bound on the MISE in the same context as above, by bounding the MISE of the projection estimator. The bias in Proposition 3 is trivially bounded within the smoothness classes defined above, so what remains to do is to bound the variance term uniformly over the same classes.

In the following we denote by $R_m$ the $m \times m$ upper-left submatrix of the infinite array $[(\Pi \mathbf{1}_k, \Pi \mathbf{1}_l)_\mathbb{H}]_{k,l \geq 0}$. Under (A1), $R_m$ is a symmetric positive definite matrix for all $m \geq 1$. For $f$ in $\mathbb{H}_+$, we denote by $\Delta_{f,m}$ the $m \times m$ diagonal matrix having entries $\pi_f \mathbf{1}_k = (f, \Pi \mathbf{1}_k)_\mathbb{H}$ on its diagonal.

THEOREM 2. *Let $f_\infty$ be in $\mathbb{H}_+$, $u = (u_m)_{m \geq 0}$ a positive decreasing sequence, $K$ and $C$ positive numbers and $r$ a nonnegative integer. Then for any positive integer $n$ and any integer $m \geq r$,*

$$(25) \quad \sup_{f \in \mathcal{C}_{f_\infty}(K,u,C,r) \cap \mathbb{H}_1} \mathbb{E}_f \|\widehat{f}_{m,n} - f\|_\mathbb{H}^2 \leq (Cu_m)^2 + \frac{K}{n} \operatorname{tr}(R_m^{-1} \Delta_{f_\infty,m}).$$

REMARK. The upper bound (25) can be optimized over all $m \geq r$. As expected, the bias term decreases and the variance bound increases as $m$ grows.

PROOF OF THEOREM 2. Pick a probability density $f$ in $\mathcal{C}_{f_\infty}(K,u,C,r)$ and depart from Proposition 3, noting that the squared bias term is bounded by $(Cu_m)^2$. Regarding the variance term, it is sufficient to consider $n = 1$. Let $\widehat{\mathbf{f}}_m$ denote the column vector of coordinates of $\widehat{f}_{m,1}$ in the basis $(\Pi \mathbf{1}_k)_{0 \leq k < m}$ of $V_m$. By Definition 1 and the definition of $R_m$, $\widehat{\mathbf{f}}_m = R_m^{-1} P_1 \mathbf{1}^{(m)}$, where $\mathbf{1}^{(m)}$ is the column vector function with entries $\mathbf{1}_k$, $0 \leq k < m$. Then

$$\begin{aligned}
\operatorname{var}_f(\widehat{f}_{m,1}) &= \mathbb{E}_f \|\widehat{f}_{m,1}\|_\mathbb{H}^2 - \|\mathbb{E}_f \widehat{f}_{m,1}\|_\mathbb{H}^2 \\
&\leq \mathbb{E}_f \widehat{\mathbf{f}}_m^T R_m \widehat{\mathbf{f}}_m - 0 \\
&= \mathbb{E}_f ((P_1 \mathbf{1}^{(m)})^T R_m^{-1} (P_1 \mathbf{1}^{(m)})) \\
&= \operatorname{tr}(R_m^{-1} \pi_f \mathbf{1}^{(m)} \mathbf{1}^{(m)T}).
\end{aligned}$$

The proof is concluded by observing that, as a positive definite symmetric matrix, $R_m^{-1}$ has positive entries on its diagonal and by noting that $\pi_f \mathbf{1}^{(m)} \mathbf{1}^{(m)T} = \Delta_{f,m} \leq \|f\|_{\infty, f_\infty} \Delta_{f_\infty,m}$. □

REMARK. Our objective here is only to provide an upper bound that is *uniform* over a given class of densities. For power series mixtures (Section 5) and mixtures of uniforms (Section 7), the bound on the variance $\operatorname{var}_f(\widehat{f}_{m,1})$ will be made more explicit by using orthogonal sequences. These bounds will then be derived directly from (12). However, they are closely related to the upper bound derived above. Indeed, let $\Phi_m$ denote the matrix $(\Phi_{k,\ell})_{0 \leq k,\ell < m}$, where $\Phi_{k,\ell}$ is as in Section 3. Observing that $(\phi_k, \phi_\ell)_{\mathbb{H}} = (\Phi_m R_m \Phi_m^T)_{k,\ell}$ for all $0 \leq k, \ell < m$, we obtain $R_m^{-1} = \Phi_m^T \Phi_m$. This relates (25) to orthonormal sequence techniques.

4.3. *Existence of smooth densities.* Theorems 1 and 2 provide lower and upper bounds, respectively, on the MISE. The classes over which these bounds apply are different in structure though; the class in Theorem 1 is a ball centered at $f_0$, while that in Theorem 2 is centered at 0. Therefore, the two bounds are not immediately comparable. The purpose of the following result is to show that under some conditions the former class is included in the latter one, thus implying that the lower bound is indeed smaller than the upper bound.

PROPOSITION 4. *Let $f_\infty$ be in $\mathbb{H}_+$, $u = (u_m)_{m \geq 0}$ a positive decreasing sequence and $r$ a nonnegative integer. Assume that we have a density $f_0$ in $\mathbb{H}_1$ and a nonnegative $C_0$ such that $f_0$ belongs to $\mathcal{C}(u, C_0, r)$.*

*Then for any positive $K$ and $K'$ satisfying $K'/(1+K) \geq \|f_0\|_{\infty,f_\infty}$ and any nonnegative $C$ and $C'$ satisfying $C' - C \geq C_0$, the inclusion*

$$(26) \qquad f_0 + \mathcal{C}_{f_0}(K, u, C, r) \subseteq \mathcal{C}_{f_\infty}(K', u, C', r)$$

*holds.*

PROOF. This follows from the inclusion $\mathcal{C}(u, C_0, r) + \mathcal{C}(u, C, r) \subseteq \mathcal{C}(u, C_0 + C, r)$ and the inequality $\|f_0 + f\|_{\infty, f_\infty} \leq (1 + \|f\|_{\infty, f_0}) \|f_0\|_{\infty, f_\infty}$. □

In the case where the inclusion (26) holds, the lower and upper bounds of Theorems 1 and 2, respectively, apply in a common setting. Hence, it is important to be able, given a smoothness class, to find $f_0$ satisfying the assumptions of Proposition 4. Under (A1), given any sequence $u = (u_m)$, it is always possible to find a nonnegative number $C_0$ such that the class $\mathcal{C}(u, C_0, r)$ contains a probability density $f_0$ for all nonnegative $r$. Take $f_0 = \Pi \mathbf{1}_0 / \nu \Pi \mathbf{1}_0$; we then trivially have $f_0 \in \mathbb{H}_1$ and $f_0 \in \mathcal{C}(u, C_0, r)$ for all $C_0 \geq 0$ if $r > 0$, or for all $C_0 \geq \|f_0\|_{\mathbb{H}}$ otherwise. This choice will indeed be made in the case of a power series mixture in Section 5. In general, this $f_0$ does not guarantee the norm $\|\cdot\|_{\infty, f_0}$ to be finite on the sets $(V_m)_{m \geq 0}$, however, which is crucial for the lower bound (see the remark following Theorem 2). In the rest of this section we provide a general construction of $f_0$ which satisfies this constraint.

Define

$$\mathbb{H}_* := \left\{ \sum_{k \geq 0} \alpha_k \Pi \mathbf{1}_k \in \mathbb{H} : \alpha_k > 0 \text{ for all } k \geq 0 \right\}.$$

By $\sum_{k\geq 0}\alpha_k\Pi\mathbf{1}_k \in \mathbb{H}$, we mean that $\sum_{k=0}^n \alpha_k\Pi\mathbf{1}_k$ converges in $\mathbb{H}$ as $n$ tends to infinity. Note that the series having nonnegative terms, by the monotone convergence theorem, this is equivalent to saying that

$$\int\left(\sum_{k\geq 0}\alpha_k\Pi\mathbf{1}_k\right)^2 d\nu < \infty.$$

Of course, $\mathbb{H}_*$ is contained in $\mathbb{H}_+$, and for any function $f = \sum_{k\geq 0}\alpha_k\Pi\mathbf{1}_k \in \mathbb{H}_*$, we have $\|\Pi\mathbf{1}_k\|_{\infty,f} \leq \alpha_k^{-1}$ for all $k$; consequently, $\|\cdot\|_{\infty,f}$ is a finite norm on every $V_m$. We now show the existence of a "smooth probability density" $f_0$ in $\mathbb{H}_*$, given any smoothness sequence $u = (u_m)$.

PROPOSITION 5. *Assume* (A1). *Then $\mathbb{H}_*$ and $\mathbb{H}_1$ have a nonempty intersection. Moreover, for any positive decreasing sequence $(u_m)$, the following holds:*

(i) *For any positive $C_0$, there are elements in $\mathbb{H}_* \cap \mathbb{H}_1$ which also belong to $\mathcal{C}(u, C_0, 1)$ and, hence, to $\mathcal{C}(u, C_0, r)$ for any positive integer $r$.*

(ii) *There exists a positive constant $C_0$ such that there are elements in $\mathbb{H}_* \cap \mathbb{H}_1$ which also belong to $\mathcal{C}(u, C_0)$.*

PROOF. The linear independence part of (A1) implies $\nu\Pi\mathbf{1}_k \neq 0$ for all $k$. For any positive sequence $(\alpha_k)$, a simple sufficient condition to have $\sum_k \alpha_k\Pi\mathbf{1}_k$ in $\mathbb{H}$ is absolute convergence, that is, $\sum_k \alpha_k\|\Pi\mathbf{1}_k\|_{\mathbb{H}} < \infty$. Moreover, by the monotone convergence theorem,

$$\nu\sum_{k\geq 0}\alpha_k\Pi\mathbf{1}_k = \sum_{k\geq 0}\alpha_k\nu\Pi\mathbf{1}_k.$$

Hence, we may pick $(\alpha_k)$ with $\alpha_k > 0$ for all $k$ and such that $\sum_k \alpha_k\Pi\mathbf{1}_k$ is both in $\mathbb{H}$ and in $L^1(\nu)$. It is then also in $\mathbb{H}_1$ by normalizing appropriately. Hence, the first part of the proposition.

For any $f = \sum_k \alpha_k\Pi\mathbf{1}_k \in \mathbb{H}_*$, since $\sum_{k=0}^{m-1}\alpha_k\Pi\mathbf{1}_k \in V_m$ and since $\mathrm{Proj}_{V_m} f$ minimizes $\|f - g\|_{\mathbb{H}}$ over $g \in V_m$, we have

$$\|f - \mathrm{Proj}_{V_m} f\|_{\mathbb{H}} \leq \left\|\sum_{k\geq m}\alpha_k\Pi\mathbf{1}_k\right\|_{\mathbb{H}} \leq \sum_{k\geq m}\alpha_k\|\Pi\mathbf{1}_k\|_{\mathbb{H}} \qquad \text{for all } m \geq 0.$$

Hence, for having $f$ in $\mathcal{C}(u, C, r) \cap \mathbb{H}_* \cap \mathbb{H}_1$, it is sufficient that $(\alpha_k)$ satisfies

(27) $\quad \sum_{k\geq 0}\alpha_k\nu\Pi\mathbf{1}_k = 1 \quad \text{and} \quad \sum_{k\geq m}\alpha_k\|\Pi\mathbf{1}_k\|_{\mathbb{H}} \leq Cu_m \qquad \text{for all } m \geq r.$

The second constraint simply says that the $\alpha_k$'s cannot be too large for $k \geq r$. If $r \geq 1$, the first constraint is then met by adapting the values of $\alpha_k$ for $k = 0, \ldots, r - 1$. If $r = 0$, then $C$ must be taken large enough for both constraints to be compatible. We now formalize these ideas.

Let $(v_m)_{m \geq 0}$ be a positive decreasing sequence such that $v_m \leq u_m$ for all $m \geq 0$ and $\lim v_m = 0$. Define a sequence $(\beta_k)$ by

$$\beta_k := (v_k - v_{k+1})(\|\Pi\mathbf{1}_k\|_{\mathbb{H}} \vee \nu\Pi\mathbf{1}_k)^{-1} \qquad \text{for all } k \geq 0.$$

Then, by construction, $(\beta_k)$ is a positive sequence and, for all $m \geq 0$, both $\sum_{k \geq m} \beta_k \|\Pi\mathbf{1}_k\|_{\mathbb{H}}$ and $\sum_{k \geq m} \beta_k \nu\Pi\mathbf{1}_k$ are less than $u_m$. Now pick a positive number $C$. Take $\alpha_k = \lambda\beta_k$ for all $k > 0$, where $0 < \lambda \leq C$ and $\lambda < (\sum_{k>0} \beta_k \nu\Pi\mathbf{1}_k)^{-1}$. Then the second part of (27) holds with $r = 1$ and we may choose $\alpha_0 > 0$ for insuring the first part of (27). It follows that $f_0 := \sum_k \alpha_k \Pi\mathbf{1}_k \in \mathbb{H}_* \cap \mathbb{H}_1 \cap \mathcal{C}(u, C, 1)$. This proves (i).

For the case $r = 0$, define $C_0 := (\sum_{k \geq 0} \beta_k \nu\Pi\mathbf{1}_k)^{-1}$; this a finite positive number by the definition of $(\beta_k)$. Putting $\alpha_k = C_0 \beta_k$ for all $k \geq 0$, (27) holds for $C \geq C_0$ and $r = 0$. This proves (ii). □

4.4. *Minimax optimality.* By optimizing the bounds (18) and (25) over $m \geq r$ in a common setting (as detailed in the previous section), we obtain lower and upper bounds on the minimax MISE over classes $\mathcal{C}_{f_\infty}(K, u, C, r)$ under the simple assumption (A1). Depending on how these bounds compare, we may obtain the minimax rate and possibly the asymptotic constant of the MISE achievable over such a class. However, this is not guaranteed. A crucial step for the lower bound is the computation of $K_{\infty, f_0}$, which will be possible only for particular smoothness classes. Concerning the upper bound, we will need to find a tractable bound on the variance, and this will only be possible in cases where orthonormal sequences are easily obtained.

In Section 5 these steps will be carried out for power series mixtures, resulting in minimax rates over smoothness classes as defined above. However, we will also give examples of mixands with different characteristics. In the setting of translation mixtures or deconvolution, treated in Section 6, an upper bound applies uniformly over all $f$ in $\mathbb{H}_1$. We will then derive a better adapted lower bound of the same rate. In the setting of mixtures of discrete uniform distributions examined in Section 7, $\Pi\mathbf{1}_k$ is not in $L^1(\nu)$ for the most natural choice of $\nu$. We will then choose $(V_m)$ different from (13) and adapt the proof of Theorem 1 to this choice. Finally, in Section 8 we give situations in which how the lower and upper bounds compare is an open question.

**5. Power series mixtures.** Let $(a_k)_{k \geq 0}$ be a sequence of positive numbers with $a_0 = 1$, and let $R$, $0 < R \leq \infty$, be the radius of convergence of the power series

$$Z(t) := \sum_{k \geq 0} a_k t^k.$$

Obviously, $Z(0) = 1$ and $Z$ is an increasing function on $[0, R)$. Put $\widetilde{Z}(t) := 1/Z(t)$. For all $\theta \in [0, R)$, the discrete distribution $\pi_\theta$ is defined by

$$(28) \qquad \pi_\theta(k) = \Pi \mathbf{1}_k(\theta) = a_k \theta^k \widetilde{Z}(\theta) \qquad \text{for all } k \geq 0.$$

In particular, the Poisson and negative binomial distributions are obtained using, respectively, $a_k = 1/k!$ and $a_k = \binom{\nu+k-1}{k}$. It is without loss of generality to assume $a_0 = 1$, since any constant multiplier of $(a_k)$ does not alter $\pi_\theta$.

Recall that $\mathbb{H} = L^2(\nu)$, where $\nu$ is a Radon measure on $\Theta$; in the case of the above power series mixture, $\Theta$ is a Borel subset of $[0, R)$. Let us first give sufficient and necessary conditions on $\nu$ for our previous results to apply, that is, for assumption (A1) to hold. These conditions are as follows.

PROPOSITION 6. *For mixands given by* (28), (A1) *is equivalent to having both the following assertions:*

(i) $\int_\Theta \theta^k \widetilde{Z}(\theta) \nu(d\theta)$ *is finite for all nonnegative integers* $k$;
(ii) $\nu$ *is not a finite sum of point masses.*

PROOF. Condition (i) exactly says that $\Pi \mathbf{1}_k$ is in $L^1(\nu)$ for all $k$. Since $\widetilde{Z}$ is bounded by one, it also gives that $\int_\Theta \theta^{2k} \widetilde{Z}^2(\theta) \nu(d\theta) < \infty$, that is, $\Pi \mathbf{1}_k$ is in $L^2(\nu)$ for all $k$. Hence, (i) is necessary and it is sufficient for having a sequence in both $L^1(\nu)$ and $L^2(\nu)$.

We now claim that the sequence $(\Pi \mathbf{1}_k)_{k \geq 0}$ is linearly independent in $\mathbb{H}$ if and only if (ii) holds. First note that if (ii) does not hold, then $\mathbb{H}$ is finite-dimensional and cannot contain an infinite sequence of linearly independent elements. To prove the converse implication, assume that (ii) holds, so that the support of $\nu$ is infinite. Pick a nonnegative integer $p$ and let $(\lambda_k)_{0 \leq k \leq p}$ be scalars such that $\sum_{0 \leq k \leq p} \lambda_k \Pi \mathbf{1}_k$ is the zero element of $\mathbb{H}$, that is, $\sum_{0 \leq k \leq p} \lambda_k \pi_\theta(k) = 0$ for $\nu$-a.e. $\theta \in \Theta$. Since $\pi.(k)$ is continuous on $\Theta$ for all $k$, $\{\theta \in \Theta : \sum_{0 \leq k \leq p} \lambda_k \pi_\theta(k) = 0\}$ is a closed set (in the relative topology on $\Theta$). Consequently, it contains the support of $\nu$ and, thus by (ii), $p+1$ distinct points $\theta_i \in \Theta$, $i = 0, \ldots, p$. As $\widetilde{Z} > 0$, it follows that $\sum_{0 \leq k \leq p} \lambda_k a_k \theta_i^k = 0$ for $i = 0, \ldots, p$, which in turn implies $\lambda_k = 0$ for $k = 0, \ldots, p$. This shows that $(\Pi \mathbf{1}_k)_{0 \leq k \leq p}$ are linearly independent for all $p \geq 0$, which completes the proof. □

The objective of the remainder of this section is to carefully apply Theorem 1 to power series mixtures when $\nu$ is Lebesgue measure on a compact interval, and to find upper bounds on the MISE for the projection estimator. This is organized as follows. We first provide computational expressions for the projection estimator in Section 5.1. We then examine the smoothness classes defined by (14), (15) and (17), and how these classes intersect $\mathbb{H}_1$ (see Section 5.2). In this context and under a submultiplicative assumption on the sequence $(a_k)$, we find that the upper and lower bounds on the MISE have the same rate, the minimax rate. A closer look is made when $R < \infty$ and also for Poisson mixtures (for which $R = \infty$). These results are stated in Section 5.3, where they are also compared to previous results found in similar settings.

5.1. *Computations based on orthonormal polynomials.* In this section we shall elaborate on the use of orthonormal polynomials in connection with the projection estimator and power series mixands. These polynomials will serve two purposes: being building blocks for numerical computations of the projection estimator and being a mathematical vehicle for establishing bounds on its variance.

The projection estimator may be computed using the techniques of Section 3. More precisely, since $V_m = \mathrm{Span}(\Pi \mathbf{1}_k, 0 \leq k < m)$ and $P_n \mathbf{1}_\ell$ is the empirical frequency of $\ell$ in the sample $(X_i)_{1 \leq i \leq n}$, (12) translates into

$$(29) \qquad \widehat{f}_{m,n} = \sum_{k=0}^{m-1} \sum_{\ell=0}^{k} \Phi_{k,\ell}(P_n \mathbf{1}_\ell)\phi_k = \frac{1}{n} \sum_{k=0}^{m-1} \sum_{i=1}^{n} \Phi_{k,X_i}\phi_k;$$

recall that $\Phi_{k,\ell} := 0$ for $\ell > k$. In the case of power series mixtures, we may use orthogonal polynomial techniques for constructing the sequence $(\phi_k)$. Let $\mathcal{P}_m$ be the set of polynomials of degree at most $m$ (with the convention $\mathcal{P}_{-1} = \{0\}$). In view of (28),

$$(30) \qquad V_m = \{p\widetilde{Z} : p \in \mathcal{P}_{m-1}\}.$$

Define the measure $\nu'$ on $\Theta$ by $d\nu' = \widetilde{Z}^2 \, d\nu$ and let $\mathbb{H}' = L^2(\nu')$. Then for any two polynomials $p$ and $q$, $(p\widetilde{Z}, q\widetilde{Z})_{\mathbb{H}} = (p,q)_{\mathbb{H}'}$. Hence, if $(q_k^{\nu'})_{k \geq 0}$ is a sequence of orthonormal polynomials in $\mathbb{H}'$ with

$$(31) \qquad q_k^{\nu'}(t) = \sum_{l=0}^{k} Q_{k,l}^{\nu'} t^l,$$

then the sequence $(\phi_k)_{k \geq 0}$ defined by

$$\phi_k(t) = q_k^{\nu'}(t)\widetilde{Z}(t) = \sum_{l=0}^{k} Q_{k,l}^{\nu'} t^l \widetilde{Z}(t) = \sum_{l=0}^{k} \frac{Q_{k,l}^{\nu'}}{a_l} \pi_t(l)$$

is an orthonormal sequence in $\mathbb{H}$ such that $(\phi_k)_{0 \leq k < m}$ spans $V_m$. Thus, $\Phi_{k,l} = (Q_{k,l}^{\nu'}/a_l)\mathbf{1}(l \leq k)$ in (29). This shows that $\widehat{f}_{m,n}$ is the same estimator as the one defined by Loh and Zhang ([17], equation (18)) with weight function $w \equiv 1$. However, it differs from the one studied by Hengartner [12], since the latter is a polynomial, and ours is in $V_m$. The coefficients $Q_{k,l}^{\nu'}$ may be obtained using standard methods to compute orthogonal sequences of polynomials; a particular method is described in the Appendix.

Let us also derive another estimator, denoted by $\check{f}_{m,n}$ and belonging to the space $\widetilde{V}_m := \{pZ : p \in \mathcal{P}_{m-1}\}$. In analogy with Definition 1, this estimator is defined as the element of $\widetilde{V}_m$ satisfying

$$(32) \qquad (\check{f}_{m,n}, g)_{\mathbb{H}} = P_n \Pi^{-1} g \qquad \text{for all } g \in V_m.$$

Observe that $(\check{f}_{m,n}, g)_{\mathbb{H}} = (\widetilde{Z}\check{f}_{m,n}, Zg)_{\mathbb{H}}$, so that (32) is equivalent to

$$(33) \qquad (\widetilde{Z}\check{f}_{m,n}, p)_{\mathbb{H}} = P_n \Pi^{-1}(\widetilde{Z}p) \qquad \text{for all } p \in \mathcal{P}_{m-1}.$$

Since $P_n \Pi^{-1}(\widetilde{Z}\cdot)$ is a linear functional on $\mathcal{P}_{m-1}$, this uniquely defines $\widetilde{Z}\check{f}_{m,n}$ in $\mathcal{P}_{m-1}$ and thus $\check{f}_{m,n}$.

We see from (32) and (7) that $\widehat{f}_{m,n} = \operatorname{Proj}_{V_m} \check{f}_{m,n}$. Therefore, by linearity, $\widehat{f}_{m,n} - \mathbb{E}_f \widehat{f}_{m,n} = \operatorname{Proj}_{V_m}(\check{f}_{m,n} - \mathbb{E}_f \check{f}_{m,n})$. Since projections do not increase the norm, taking the squared norm and expectation gives

$$(34) \qquad \operatorname{var}_f(\widehat{f}_{m,n}) \leq \operatorname{var}_f(\check{f}_{m,n}) \qquad \text{for all } f \in \mathbb{H}_1.$$

We will use this property below to bound the variance of $\widehat{f}_{m,n}$. At the moment let us note that this bound indicates that $\check{f}_{m,n}$ does not behave as well as $\widehat{f}_{m,n}$, even though, for brevity, we leave aside the problem of the bias. Nevertheless, the estimator $\check{f}_{m,n}$ has the appealing property that it may be expressed by using a sequence $(q_k^\nu)_{k \geq 0}$ of orthonormal polynomials in $\mathbb{H} = L^2(\nu)$ that does not depend on the sequence $(a_k)$ but only on $\nu$. To see this, let us write, as in (31),

$$q_k^\nu(t) = \sum_{l=0}^{k} Q_{k,l}^\nu t^l.$$

Again, by convention, we extend the values of $Q_{k,l}^\nu$ to the domain $l > k$ by zeros. An algorithm for computing $Q_{k,l}^\nu$ is given in the Appendix for $\Theta = [a,b]$ and certain choices of $\nu$. Let us now express $\check{f}_{m,n}$ in terms of this sequence. By (33), as $\Pi^{-1}(\widetilde{Z}(t)t^l) = \Pi^{-1}(\pi.(l)/a_l) = \mathbf{1}_l/a_l$, we obtain

$$(\widetilde{Z}\check{f}_{m,n}, q_k^\nu)_{\mathbb{H}} = \sum_{\ell=0}^{k} \frac{Q_{k,\ell}^\nu}{a_\ell} P_n \mathbf{1}_\ell.$$

Since $\widetilde{Z}\check{f}_{m,n}$ belongs to $\mathcal{P}_{m-1}$, we conclude that the right-hand side of this display has the coefficients of the expansion of $\widetilde{Z}\check{f}_{m,n}$ in the orthonormal basis $(q_k^\nu)$. Thus,

$$(35) \qquad \check{f}_{m,n} = Z \sum_{k=0}^{m-1} \sum_{\ell=0}^{k} \frac{Q_{k,\ell}^\nu}{a_k}(P_n \mathbf{1}_\ell) q_k^\nu = \frac{Z}{n} \sum_{k=0}^{m-1} \sum_{i=1}^{n} \frac{Q_{k,X_i}^\nu}{a_{X_i}} q_k^\nu.$$

Below we will use this expression to bound the variance of $\widehat{f}_{m,n}$.

5.2. *Approximation classes.* Recall the definitions (14), (15) and (17) made in Section 4. These smoothness classes are closely related to those used in previous works on power series mixtures, as will be shown in this section. The discussion will be devoted to the case

$\nu$ is Lebesgue measure on $[a,b]$ with $0 \leq a < b < R$.

Hence, $\mathbb{H}$ is the usual $L^2$ space of functions on $[a,b]$. Let, for any positive $\alpha$,

$$\mathbf{u}^\alpha := ((1+n)^{-\alpha})_{n \geq 0}.$$

It turns out that, for these particular sequences, the classes $\mathcal{C}(\mathbf{u}^\alpha, C)$ of Section 4 are equivalent to classes defined using weighted moduli of smoothness. This, in turn, will relate them to Sobolev and Hölder classes, classes that were considered by Hengartner [12]. To make this precise, let $\|\cdot\|_p$ be the $L^p$ norm over $[a,b]$, define the function $\phi(x) = \sqrt{(x-a)(b-x)}$ on this interval and let $\Delta_h^r(f,x)$ be the symmetric difference of order $r$, that is,

$$\Delta_h^r(f,x) := \sum_{i=0}^r \binom{r}{i} (-1)^i f(x + (i-r/2)h),$$

with the classical convention that $\Delta_h^r(f,x)$ is set to 0 if $x+(i-r/2)h$ is outside $[a,b]$ for $i=0$ or $r$. Then for any function $f$ on $[a,b]$, the weighted modulus of smoothness is defined as

(36) $$\omega_r^\phi(f,t)_p := \sup_{0<h\leq t} \|\Delta_{h\phi(\cdot)}^r(f,\cdot)\|_p.$$

The effect of the weight $\phi$ here is to relax the regularity conditions on $f$ at the endpoints $a$ and $b$. Finally, for all positive numbers $\alpha$ and $C$, define the classes

(37) $$\widetilde{\mathcal{C}}(\alpha, C) := \{f \in \mathbb{H} : \|f\|_\mathbb{H} \leq C, \omega_{[\alpha]+1}^\phi(f,t)_2 \leq Ct^\alpha \text{ for all } t > 0\}.$$

The following result shows that these classes are, in a certain sense, equivalent to $\mathcal{C}(\mathbf{u}^\alpha, C)$.

PROPOSITION 7. *For any positive number $\alpha$, there exist positive constants $C_1$ and $C_2$ such that*

(38) $$\mathcal{C}(\mathbf{u}^\alpha, C_1 C) \subseteq \widetilde{\mathcal{C}}(\alpha, C) \subseteq \mathcal{C}(\mathbf{u}^\alpha, C_2 C) \qquad \text{for all } C > 0.$$

Before giving the proof of this proposition, we explain the point of this result and of defining the classes $\widetilde{\mathcal{C}}(\alpha, C)$. Recall that the standard modulus of smoothness

$$\omega_r(f,t)_p := \sup_{0<h\leq t} \|\Delta_h^r(f,\cdot)\|_p$$

provides definitions of semi-norms for Besov and Sobolev spaces (see, resp. equation (2.10.1) and Theorem 2.9.3 in [4]). In particular, the classes defined as in (37) but with the weighted modulus of smoothness (36) replaced by the standard one with the same parameters $r = [\alpha] + 1$ and $p = 2$ are balls in the Besov space $B^\alpha_\infty(L^2[a,b])$. Now using Theorems 6.2.4 and 6.6.2 in [4] and the fact that $\|f\|_{W^r_p(\phi)} \leq \|f\|_{W^r_p}$, from [4], equation (6.6.5), we have that, for a constant $C_0 > 0$ only depending on $p$ and $r$,

$$\omega^\phi_r(f,t)_p \leq C_0 \omega_r(f,t)_p \quad \text{for } 0 < t < (2r)^{-1}.$$

Furthermore, bounding $\omega^\phi_r(f,t)_p$ by $\|f\|_p$ up to a multiplicative constant for $t \geq (2r)^{-1}$ as in [4], equation (6.6.5) shows that $\widetilde{\mathcal{C}}(\alpha, C)$ contains Besov balls $\{\|f\|_{B^\alpha_\infty(L^2[a,b])} \leq C'_0 C\}$ for a constant $C'_0$, but is not contained in such balls. Using inequalities between Hölder, Sobolev and Besov semi-norms, it also follows that $\widetilde{\mathcal{C}}(\alpha, C)$ contains balls of the Hölder space $C^\alpha[a,b]$ and of the Sobolev space $W^\alpha_2$ and, of course, converse inclusions are not to be found. In view of Proposition 7, since $\widetilde{\mathcal{C}}(\alpha, C)$ contains Besov, Hölder and Sobolev balls as just described, so does $\mathcal{C}(\mathbf{u}^\alpha, C)$.

These inclusions are helpful for comparing our results to those of Hengartner [12], where minimax rates are given for Sobolev balls with integer exponents and conjectured for Hölder balls. In his paper, as well as the present one, rates for the projection estimator are obtained using properties which hold over the classes $\widetilde{\mathcal{C}}(\alpha, C)$ and, consequently, over smaller ones such as Sobolev and Hölder balls. Minimax bounds, however, are obtained using different methods. Our approach takes advantage of the whole class over which the rate applies, whereas Hengartner [12] only used subclasses to derive minimax bounds. This "closer look" allows us to derive minimax bounds applying to $\widetilde{\mathcal{C}}(\alpha, C)$ for all $\alpha \geq 1$, not only integers, and to obtain results on the asymptotic constant when refining the class $\widetilde{\mathcal{C}}(\alpha, C)$ to $\mathcal{C}(u, C, r)$.

PROOF OF PROPOSITION 7. Write the equivalence relationships (38) as

$$\mathcal{C}(\mathbf{u}^\alpha, \cdot) \asymp \widetilde{\mathcal{C}}(\alpha, \cdot).$$

We start by relating $\widetilde{\mathcal{C}}(\alpha, C)$ to classes of the form

$$\overline{\mathcal{C}}(u, C) := \left\{ f \in \mathbb{H} : \inf_{p \in \mathcal{P}_{m-1}} \|f - p\|_\mathbb{H} \leq C u_m \text{ for all } m \geq 0 \right\};$$

recall that $\mathcal{P}_m$ is the set of polynomials of degree at most $m$. Theorem 8.7.3 and equation (8.7.25) of [4] show that, for all $\alpha > 0$ and all $r > \alpha$, there exist constants $C'_1$ and $C'_2$ such that, for all $C > 0$ and $f \in \mathbb{H}$,

$$\sup_{t \geq r} \omega^\phi_r(f, 1/t)_2 t^\alpha \leq C \implies \sup_{m \geq 4r} \inf_{p \in \mathcal{P}_m} \|f - p\|_\mathbb{H} m^\alpha \leq C'_1 C,$$

$$\sup_{m \geq r} \inf_{p \in \mathcal{P}_m} \|f - p\|_\mathbb{H} m^\alpha \leq C'_2 C \implies \sup_{t > r} \omega^\phi_r(f, 1/t)_2 t^\alpha \leq C.$$

Here $C_1'$ and $C_2'$ may depend on $\alpha$ and $r$. Taking $r = [\alpha] + 1$, observing that $(1+m)^\alpha \asymp (m-1)^\alpha$ for $m \geq 2$ and using $\|f\|_{\mathbb{H}}$ for bounding $\inf_{p \in \mathcal{P}_{m-1}} \|f - p\|_{\mathbb{H}}$ and $\omega_r^\phi(f, 1/t)_2$ in cases not covered by the above implications, we obtain $\widetilde{\mathcal{C}}(\alpha, \cdot) \asymp \overline{\mathcal{C}}(\mathbf{u}^\alpha, \cdot)$. Thus, also $Z\widetilde{\mathcal{C}}(\alpha, \cdot) \asymp Z\overline{\mathcal{C}}(\mathbf{u}^\alpha, \cdot)$, where $Z\widetilde{\mathcal{C}}(\alpha, C) := \{Zf : f \in \widetilde{\mathcal{C}}(\alpha, C)\}$ and so on.

Next we proceed to study $Z\overline{\mathcal{C}}(\mathbf{u}^\alpha, C)$. By (14) and (30),

$$\mathcal{C}(u, C) = \left\{ f \in \mathbb{H} : \inf_{p \in \mathcal{P}_{m-1}} \|f - p\widetilde{Z}\|_{\mathbb{H}} \leq Cu_m \text{ for all } m \geq 0 \right\}.$$

Since $\widetilde{Z}$ is positive and decreasing on $[a, b]$,

(39) $\qquad \widetilde{Z}(b)\|f\|_{\mathbb{H}} \leq \|\widetilde{Z}f\|_{\mathbb{H}} \leq \widetilde{Z}(a)\|f\|_{\mathbb{H}} \qquad$ for all $f$ in $\mathbb{H}$.

This shows that $\|f - p\widetilde{Z}\|_{\mathbb{H}} \asymp \|Zf - p\|_{\mathbb{H}}$, whence $\mathcal{C}(\mathbf{u}^\alpha, \cdot) \asymp Z\overline{\mathcal{C}}(\mathbf{u}^\alpha, \cdot)$. Recalling that $Z\widetilde{\mathcal{C}}(\alpha, \cdot) \asymp Z\overline{\mathcal{C}}(\mathbf{u}^\alpha, \cdot)$, we thus see that in order to prove (38) it is sufficient to show $Z\widetilde{\mathcal{C}}(\alpha, \cdot) \asymp \widetilde{\mathcal{C}}(\alpha, \cdot)$.

The remainder of the proof is thus devoted to showing that there are constants $C_1'$ and $C_2'$ such that $f \in \widetilde{\mathcal{C}}(\alpha, C)$ implies $Zf \in \widetilde{\mathcal{C}}(\alpha, C_1'C)$ and $\widetilde{Z}f \in \widetilde{\mathcal{C}}(\alpha, C_2'C)$. Since $b < R$, $Z$ is bounded away from zero and infinity on $[a, b]$ and both $Z$ and $\widetilde{Z}$ are thus infinitely continuously differentiable on this interval. Having made this observation, both of the desired implications follow from the claim that, for any $[\alpha] + 1$ times continuously differentiable function $g$ on $[a, b]$, there exists $c > 0$ such that

(40) $\qquad f \in \widetilde{\mathcal{C}}(\alpha, C) \implies gf \in \widetilde{\mathcal{C}}(\alpha, cC).$

To prove this claim, pick an $f$ in $\widetilde{\mathcal{C}}(\alpha, C)$ and let $r := [\alpha] + 1$. Recalling that $\widetilde{\mathcal{C}}(\alpha, \cdot) \asymp \overline{\mathcal{C}}(\mathbf{u}^\alpha, \cdot)$, we see that the union $\bigcup_{c'>0} \widetilde{\mathcal{C}}(\alpha, c')$ coincides with $\bigcup_{c'>0} \overline{\mathcal{C}}(\mathbf{u}^\alpha, c')$, and is, hence, increasing as $\alpha$ decreases. As the union can be written $\bigcup_{c'>0} \widetilde{\mathcal{C}}(\alpha, c'C)$, there exists a positive $c'$, depending only on $\alpha$, such that $\widetilde{\mathcal{C}}(\alpha, C) \subseteq \widetilde{\mathcal{C}}(\alpha - r + i, c'C)$ for all $i = 1, \ldots, r$. Since $f$ is included in all these classes and $r = [\alpha] + 1$, we find that $\omega_i^\phi(f, t)_2 = \omega_{[\alpha-r+i]+1}^\phi(f, t)_2 \leq c'Ct^{\alpha-r+i}$ for these $i$, and also for $i = 0$ with the usual convention $\omega_0^\phi(f, t)_2 := \|f\|_{\mathbb{H}}$.

Now the equality (obtained by standard algebra)

$$\Delta_h^r(fg, x) = \sum_{i=0}^r \binom{r}{i} \Delta_h^i(f, x + (r-i)h/2) \Delta_h^{r-i}(g, x - ih/2)$$

and the bound $|\Delta_h^{r-i}(g, x - ih/2)| \leq \|g^{(r-i)}\|_{L^\infty[a,b]} (rh)^{r-i}$ for all $0 \leq i \leq r$ and $x \in [a, b]$ yield

$$\omega_r^\phi(fg, t)_2 \leq M_g \sum_{i=0}^r \binom{r}{i} \omega_i^\phi(f, t)_2 (rt)^{r-i},$$

where $M_g := \max_{0 \leq j \leq r} \|g^{(j)}\|_{L^\infty(a,b)}$. Since, as shown above, $\omega_i^\phi(f,t)_2 t^{r-i} \leq c' C t^\alpha$, the claim (40) follows with $c = c'(1+r)^r M_g$. □

Let us now consider the smoothness classes $\mathcal{C}_{f_0}(K, u, C, r)$ defined in (17). We take $f_0$ such that, for two positive constants $c_1$ and $c_2$,

(41) $$c_1 \leq f_0(t) \leq c_2 \quad \text{for all } t \in [a,b].$$

Under this condition the norm (16) satisfies

(42) $$\frac{1}{c_2} \operatorname*{ess\,sup}_{a \leq t \leq b} |f(t)| \leq \|f\|_{\infty, f_0} \leq \frac{1}{c_1} \operatorname*{ess\,sup}_{a \leq t \leq b} |f(t)|.$$

Thus, as in (38), the classes defined by (17) are equivalent to classes defined by the weighted modulus of smoothness and a bound on the sup norm. Indeed, with $\widetilde{\mathcal{C}}(K, \alpha, C) := \{f \in \widetilde{\mathcal{C}}(\alpha, C) : \operatorname{ess\,sup}_{a \leq t \leq b} |f(t)| \leq K\}$,

(43) $$\mathcal{C}_{f_0}(K/c_2, \mathbf{u}^\alpha, C_1 C) \subseteq \widetilde{\mathcal{C}}(K, \alpha, C) \subseteq \mathcal{C}_{f_0}(K/c_1, \mathbf{u}^\alpha, C_2 C).$$

Another important consequence of (42) is that for bounding $K_{\infty, f_0}(V_m)$ we may use the Nikolskii inequality (see, e.g., [4], Theorem 4.2.6). This inequality states that there is a universal positive constant $C$ such that, for all nonnegative integers $m$,

$$\sup\left\{\sup_{-1 \leq t \leq 1} |p(t)| : p \in \mathcal{P}_{m-1}, \int_{-1}^{1} |p(t)|^2 \, dt = 1\right\} \leq Cm.$$

Also recall (30), that is, $V_m = \{p\widetilde{Z} : p \in \mathcal{P}_{m-1}\}$. Combining these observations with (39) and the Nikolskii inequality yields, for all $m \geq 1$,

(44) $$K_{\infty, f_0}(V_m) \leq C_{a,b} m$$

for a positive constant $C_{a,b}$ depending only on $(a_k)$, $a$, $b$ and $c_2$.

The following useful result says how the class $\widetilde{\mathcal{C}}(\alpha, C)$ intersects $\mathbb{H}_1$.

LEMMA 3. *Let $\alpha$ be a positive number. If $C < 1/\sqrt{b-a}$, then the intersection of $\widetilde{\mathcal{C}}(\alpha, C)$ with $\mathbb{H}_1$ is empty. Furthermore, the intersection of $\widetilde{\mathcal{C}}(\alpha, 1/\sqrt{b-a})$ with $\mathbb{H}_1$ is the singleton set $\{\mathbf{1}_{[a,b]}/(b-a)\}$.*

PROOF. Pick $f$ in $\mathbb{H}_1$. Applying Jensen's inequality $\mathbb{E}g(Y) \geq g(\mathbb{E}(Y))$ with $g(t) = t^2$, $Y = f$ and probability measure $dt/(b-a)$ on $[a,b]$ gives

$$\left\|\frac{f}{\sqrt{b-a}}\right\|_{\mathbb{H}} \geq \int_a^b \frac{f(t)}{b-a} \, dt = \frac{1}{b-a},$$

so that $\|f\|_{\mathbb{H}} \geq 1/\sqrt{b-a}$. Hence, the first part of the lemma. Now, using the strict convexity of the square function, equality in the above relation implies that $f$ is constant. Thus, to prove the second part of the lemma, we need to check that the uniform density $\mathbf{1}_{[a,b]}/(b-a)$ belongs to $\widetilde{\mathcal{C}}(\alpha, 1/\sqrt{b-a})$ for all $\alpha > 0$. This is trivially true since $\omega_k^\phi(\mathbf{1}_{[a,b]}, t)_2 = 0$ for all $t > 0$ and $k > 0$.
□

We conclude this section with a remark on the somewhat more general case when $\Theta = [a,b]$ and $\nu(dt) = dt/w(t)$ for a weight function $w$, investigated by Loh and Zhang [17]. In this case the classes denoted by $\mathcal{G}(\alpha, m, M, w_0)$ in [17] are included in the classes $\mathcal{C}(\mathbf{u}^\alpha, C, m)$ as

$$\mathcal{G}(\alpha, m, M, w_0) \subseteq \mathcal{C}(\mathbf{u}^\alpha, M_1 M, m+1) \qquad \text{for all } \alpha > 0,$$

$$\mathcal{G}(\alpha', m, M, w_0) \supseteq \mathcal{C}(\mathbf{u}^\alpha, M_2 M, m+1) \qquad \text{for all } \alpha > \alpha' > 0,$$

for positive constants $M_1$ and $M_2$ depending on $m$, $\alpha$ and $\alpha'$. Hence, our setting is very close to the one adopted by Loh and Zhang [17] in their Section 3, where they provide lower and upper bounds on the MISE over these classes. However, all their results in this section rely on special conditions, namely, their (19) and (20), which imply restrictions on the parameters $\alpha$ and $m$ defining the classes $\mathcal{G}(\alpha, m, M, w_0)$ (see their Remark 3). In particular, the rate optimality in these classes is only obtained for integer $\alpha$ (we refer to the closing comments of Section 3 in [17]). As remarked above, a similar restriction applies in [12], where minimax rates are proved in Sobolev classes with integer exponents. In contrast, the lower bound of Theorem 1 will provide the minimax rate for all $\alpha \geq 1$ in our classes and we will also obtain results on the asymptotic constant.

5.3. *Minimax MISE rates.* The following result is concerned with the asymptotic properties of the projection estimator and lower bounds on the MISE over the approximation classes of Section 5.2.

THEOREM 3. *Assume that $\nu$ is Lebesgue measure on $[a,b]$ with $0 \leq a < b < R$, and let $\lambda := \gamma + \sqrt{\gamma^2 + 1}$ with $\gamma = (2+a+b)/(b-a)$. Then the following assertions hold true:*

(a) *Let $\alpha$ and $C$ be positive numbers, $r$ be a nonnegative integer and $(m_n)$ be a nondecreasing divergent integer sequence. If there exists a number $\lambda_1$ larger than $\lambda$ such that*

(45) $$\frac{1}{n} \lambda_1^{2m_n} \max_{0 \leq k < m_n} \frac{b^k}{a_k} \to 0 \qquad \text{as } n \to \infty,$$

*then*

$$\sup_{f \in \mathcal{C}(\mathbf{u}^\alpha, C, r) \cap \mathbb{H}_1} \mathbb{E}_f \|\widehat{f}_{m_n, n} - f\|_{\mathbb{H}}^2 \leq C^2 m_n^{-2\alpha}(1 + o(1)).$$

(b) *Let $\alpha \geq 1$, $C$ be a positive number, $r$ be a positive integer and $(m'_n)$ be a nondecreasing divergent integer sequence. Put*

$$(46) \qquad w_n := n \sum_{k \geq m'_n} \nu \Pi \mathbf{1}_k \qquad \text{for any positive integer } n.$$

*If $(w_n)$ tends to zero, then*

$$(47) \qquad \inf_{\hat{f} \in \mathcal{S}_n} \sup_{f \in \mathcal{C}(\mathbf{u}^\alpha, C, r) \cap \mathbb{H}_1} \mathbb{E}_f \|\hat{f} - f\|_\mathbb{H}^2 \geq C^2 m'_n{}^{-2\alpha}(1 + o(1)).$$

(c) *Let $\alpha \geq 1$ and $C > 1/\sqrt{b-a}$. If there exist sequences $(m_n)$ and $(m'_n)$ satisfying the conditions of* (a) *and* (b) *and such that $\liminf_{n \to \infty} m_n/m'_n > 0$, then the minimax MISE rate over $\widetilde{\mathcal{C}}(\alpha, C) \cap \mathbb{H}_1$ is $m_n^{-2\alpha}$ and it is achieved by the projection estimator $\widehat{f}_{m_n, n}$.*

Before giving the proof in Section 5.4, we make the following remarks and examine the particular cases of mixands with $R < \infty$ and Poisson mixtures:

(i) The condition on $C$ in (c) is necessary since otherwise the class is empty or reduces to one element; see Lemma 3. In contrast, under the assumptions of (c), a direct application of (a) and (b) provides the same minimax rate over $\mathcal{C}(\mathbf{u}^\alpha, C, r) \cap \mathbb{H}_1$ for $\alpha > 1$, $r \geq 1$ and $C > 0$.

(ii) The same lower and upper bounds apply to classes defined by adding a bound on the uniform norm. For instance, part (c) holds when replacing $\widetilde{\mathcal{C}}(\alpha, C)$ by $\widetilde{\mathcal{C}}(K, \alpha, C)$ for any $K > 1$ (for $K \leq 1$, this class is empty or reduces to the uniform density as for a too small $C$). This can be verified easily by reading the proof.

(iii) The $o$-terms in parts (a) and (b) can be made more precise. In part (a), $m_n^{-2\alpha}(1 + o(1))$ can be replaced by $(m_n + 1)^{-2\alpha} + \kappa \xi^{-m_n}$, where $\xi > 1$ and $\kappa > 0$ depend only on $(a_k)$, $a$ and $b$ and the inequality holds for $m_n \geq r$. In part (b), $m'_n{}^{-2\alpha}(1 + o(1))$ can be replaced by $(m'_n + 2)^{-2\alpha} \exp(-\kappa' w_n)$, where $\kappa' > 0$ depends only on $(a_k)$, $a$ and $b$, and the inequality holds for $n$ sufficiently large.

(iv) The estimator $\widehat{f}_{m_n, n}$ in (a) depends only on $(m_n)$, which is fixed by $(a_k)$, $a$ and $b$ through (45). Thus, it is *universal* in the sense that it does not depend on $C$ or $\alpha$, although, under the conditions of (c), it is rate optimal in $\widetilde{\mathcal{C}}(\alpha, C) \cap \mathbb{H}_1$ for all $\alpha \geq 1$ and all $C > 1/\sqrt{b-a}$. However, an interesting problem, which is left open in the present work, would be to build an estimator which adapts to an unknown $[a, b] \subset [0, R)$.

(v) Clearly, one can always find two sequences $(m_n)$ and $(m'_n)$ satisfying the conditions of (a) and (b), respectively. In contrast, for obtaining (c), it remains to show that these sequences can be chosen equivalent. This requires further conditions on the asymptotics of $(a_k)$.

A condition addressing the issue of item (v) above is the following:

(A2)    There exists a positive constant $c_0$ such that $a_{k+l} \leq c_0 a_k a_l$ for all nonnegative integers $k$ and $l$.

This condition holds in various important situations including Poisson and negative binomial mixands. In Proposition 8 below we show that it ensures that Theorem 3(c) applies. The condition says that $(a_k)$ is submultiplicative up to the constant $c_0$. Mimicking the argument of the subadditive lemma (see, e.g., [3], page 231, and Exercise 6.1.17, page 235), $L = \lim_{n\to\infty} n^{-1} \log a_n$ exists and is given by $L = \inf_{n\geq 1} n^{-1}(\log c_0 + \log a_n)$. Thus, $c_0 a_n \geq e^{Ln}$ for all $n \geq 1$. Note that $L = -\infty$ is possible. Since $a_n = O(e^{n(L+\varepsilon)})$ and $e^{nL} = O(a_n)$ for all positive $\varepsilon$, $L$ is related to the radius of convergence through the relation $Re^L = 1$, that is, $L = -\infty$ if and only if $R = \infty$ and $L = -\log R$ otherwise. In addition, for $R < \infty$, we see that the series $\sum a_k \theta^k$ is divergent at $\theta = R$.

A first simple application of this assumption is the following lemma:

LEMMA 4. *Under* (A2), *we have, for any nonnegative integer* $m$,
$$\sum_{k \geq m} \nu \Pi \mathbf{1}_k \leq c_0 a_m \int t^m \nu(dt).$$

PROOF. A direct application of (A2) with (28) shows that, for all $k \geq m$, $\Pi \mathbf{1}_k(t) \leq c_0 a_m t^m \Pi \mathbf{1}_{k-m}(t)$. Thus, by monotone convergence and the observation $\sum_{k \geq 0} \Pi \mathbf{1}_k = 1$,
$$\sum_{k \geq m} \nu \Pi \mathbf{1}_k = \nu \sum_{k \geq m} \Pi \mathbf{1}_k \leq c_0 a_m \int t^m \sum_{k \geq 0} \Pi \mathbf{1}_k(t) \nu(dt) = c_0 a_m \int t^m \nu(dt). \quad \square$$

PROPOSITION 8. *Under* (A2) *there exist a sequence* $(m_n)_{n\geq 0}$ *satisfying the conditions of Theorem* 3(a) *and a number* $\eta \geq 1$ *such that, by setting* $m'_n := [\eta m_n]$, *the sequence* $(m'_n)_{n\geq 0}$ *satisfies the conditions of Theorem* 3(b). *Hence, Theorem* 3(c) *applies.*

We note that, with $m'_n := [\eta m_n]$, the asymptotic constants of parts (a) and (b) of Theorem 3 differ by a factor $\eta^{2\alpha}$. Thus, up to this factor, the projection estimator is minimax MISE efficient over classes, $\mathcal{C}(\mathbf{u}^\alpha, C, r) \cap \mathbb{H}_1$ with $\alpha > 1$ and $r \geq 1$. How large $\eta$ needs to be taken depends on the model through $(a_k)$, $a$ and $b$. The following result is a sharpened version of Proposition 8 in the case where $R < \infty$, obtained by optimizing with respect to $\eta$. We refer to the proof for details. The proof of Proposition 8 in the case where $R = \infty$ is postponed to Section 5.5. It provides an explicit, although more involved, construction of $(m_n)$ and $\eta$ in a way making $\eta \geq 2$ necessary.

COROLLARY 1. *Let $\nu$ be Lebesgue measure on $[a,b]$ with $0 \leq a < b < R$. Assume that $R$ is finite and that (A2) holds. Let $\alpha \geq 1$ and $C > 1/(b-a)$. Then the minimax MISE rate over $\widetilde{\mathcal{C}}(\alpha, C) \cap \mathbb{H}_1$ is $(\log n)^{-2\alpha}$. This rate is achieved by the projection estimator $\widehat{f}_{m_n}$ with $m_n := [\tau \log n]$ for any positive $\tau$ less than $1/\log\{\lambda^2 (1 \vee bR)\}$. Moreover, if $\alpha > 1$, then for any positive $C$ and any positive integer $r$,*

$$(48) \quad \limsup_{n\to\infty}(\log n)^\alpha \sup_{f \in \mathcal{C}(\mathbf{u}^\alpha, C, r) \cap \mathbb{H}_1} \mathbb{E}_f \|\widehat{f}_{m_n,n} - f\|_\mathbb{H} \leq \tau^{-\alpha} C$$

*and*

$$(49) \quad \liminf_{n\to\infty}(\log n)^\alpha \inf_{\hat f \in \mathcal{S}_n} \sup_{f \in \mathcal{C}(\mathbf{u}^\alpha, C, r) \cap \mathbb{H}_1} \mathbb{E}_f \|\hat f - f\|_\mathbb{H} \geq (\log(R/b))^\alpha C.$$

PROOF. Put $\tau_{\max} := 1/\log\{\lambda^2(1 \vee bR)\}$ and consider (45). We will make use of the properties derived above from (A2). First we note that, since $a_k \geq c_0^{-1} e^{Lk} = c_0^{-1} R^{-k}$ for all $k$, we have

$$\frac{1}{n}\lambda_1^{2m} \max_{0 \leq k < m} \frac{b^k}{a_k} \leq \frac{1}{n}\lambda_1^{2m} c_0 \max_{0 \leq k < m} (bR)^k = \frac{1}{n}\lambda_1^{2m} c_0 (1 \vee bR)^{m-1}.$$

Thus, the log of the left-hand side of this equation is at most $m\log(\lambda_1^2(1 \vee bR)) - \log n + \log c_0$, so that for $m_n = [\tau \log n]$ with $\tau \log(\lambda_1^2(1 \vee bR)) < 1$, (45) holds. Combining this condition with the requirement $\lambda_1 > \lambda$, we obtain the bound $\tau < \tau_{\max}$. Hence, for such $(m_n)$, Theorem 3(a) applies and gives (48).

Now consider part (b) of Theorem 3. Lemma 4 shows that $\sum_{k \geq m} \nu \Pi \mathbf{1}_k = O(a_m b^m)$. Moreover, for any $\varepsilon > 0$, it holds that $a_m \leq (e^L/(1-\varepsilon))^m = (R(1-\varepsilon))^{-m}$ for large $m$. Thus, with $m'_n = [\eta m_n]$ and $(w_n)$ defined as in (46), $\log w_n \leq \log n + \eta m_n \log(b/(R-\varepsilon))$ up to an additive constant. Hence, we may choose $\eta > -1/(\tau \log(b/(R-\varepsilon))) > 0$ such that $w_n \to 0$. This achieves the proof of the minimax MISE rate by applying Theorem 3(c) with the chosen sequences $(m_n)$ and $(m'_n)$. Theorem 3(b) gives a lower bound on the MISE asymptotically equivalent to $C^2(\eta\tau \log n)^{-2\alpha}$. Optimizing with respect to $\eta\tau$ under the above constraints and letting $\varepsilon$ tend to zero gives (49). □

We note that Loh and Zhang [17] proved the rate $(\log n)^\alpha$ to be minimax over different smoothness classes, but also under different assumptions on $(a_k)$. Corollary 1 extends their results to other classes and mixands, but only for $b < R < \infty$.

We now consider the Poisson case. We already know from Proposition 8 that we can find a universal projection estimator whose rate is optimal in classes $\widetilde{\mathcal{C}}(\alpha, C)$, which complements the results of Hengartner [12] and Loh and Zhang [17]. This does not address asymptotic efficiency, however, which includes computations of asymptotic constants. It turns out that a direct computation of $(m_n)$ provides the following precise result for $\mathcal{C}(\mathbf{u}^\alpha, C, r)$.

COROLLARY 2. *Let $\nu$ be Lebesgue measure on $[a,b]$ with $0 \leq a < b < R$ and suppose that $a_k = 1/k!$ (Poisson mixands, $R = \infty$). Let $\alpha \geq 1$ and $C > 1/\sqrt{b-a}$. Then the minimax MISE rate over $\widetilde{\mathcal{C}}(\alpha, C) \cap \mathbb{H}_1$ is $(\log n / \log\log n)^{-2\alpha}$. This rate is achieved by the projection estimator $\widehat{f}_{m_n}$ with $m_n := [\tau \log n / \log\log n]$ for any positive $\tau \leq 1$. Moreover, if $\alpha > 1$, then for any positive $C$ and any positive integer $r$, the projection estimator $\widehat{f}_{m_n}$ defined as above with $\tau = 1$ is asymptotically minimax efficient (including the constant) over $\mathcal{C}(\mathbf{u}^\alpha, C, r)$:*

$$\limsup_{n\to\infty} \left(\frac{\log n}{\log\log n}\right)^\alpha \sup_{f \in \mathcal{C}(\mathbf{u}^\alpha, C, r) \cap \mathbb{H}_1} \mathbb{E}_f \|\widehat{f}_{m_n, n} - f\|_\mathbb{H}$$
$$= \liminf_{n\to\infty} \left(\frac{\log n}{\log\log n}\right)^\alpha \inf_{\hat{f} \in \mathcal{S}_n} \sup_{f \in \mathcal{C}(\mathbf{u}^\alpha, C, r) \cap \mathbb{H}_1} \mathbb{E}_f \|\hat{f} - f\|_\mathbb{H}$$
$$= C.$$

PROOF. Consider (45). By Stirling's formula, $\max_{0 \leq k < m}(b^k/a_k) = O(m^m c^m)$ for a positive $c$. Since $m_n \leq \tau \log n / \log\log n$ for $\tau$ in $(0, 1]$, a simple computation yields

$$\log(m_n^{m_n} c^{m_n}/n) \leq (\tau - 1)\log n - (\tau + o(1))\frac{\log n \log\log\log n}{\log\log n} \to -\infty.$$

Condition (45) follows for any $\lambda_1 > \lambda$ (indeed, $\lambda_1$ simply multiplies $c$).

Now assume $\tau = 1$ and consider part (b) of Theorem 3. Use Lemma 4 once again to see that $\sum_{k \geq m} \nu \Pi \mathbf{1}_k = O(a_m b^m)$. In the present case $a_m b^m = O(m^{-m} c^m)$ for a positive $c$ (not the same as above). For any $\sigma$ in $(0, 1)$, it holds that $m_n \geq \sigma \log n / \log\log n$ for large $n$. As usual, we set $m'_n := [\eta m_n]$ for a positive number $\eta$ to be optimized later on and define $(w_n)$ as in (46). Thus, if $\eta \sigma > 1$, then up to an additive constant,

$$\log w_n \leq \log(n(\eta m_n)^{-\eta m_n} c^{\eta m_n}) \leq (1 - \eta\sigma + o(1))\log n \to -\infty.$$

The conclusions of the corollary now follow by applying the various parts of Theorem 3. In particular, the lower bound on the asymptotic MISE, normalized by the rate $(\log n / \log\log n)^{2\alpha}$, is obtained upon observing that we may choose $\sigma$ and, hence, also $\eta$, arbitrarily close to 1. $\square$

We recall that, in contrast to Corollary 2, Hengartner [12] did not consider constants and the exact rate was proved for Sobolev classes with entire exponents only. Likewise, Theorem 5 in [17] does not provide constants, and optimal rates are obtained only in cases similar to Hengartner [12] (cf. the paragraph ending Section 5.2 above). By determining the asymptotic constant, we also answer a question raised by Hengartner ([12], page 921, Remark 4). He suggested that the optimal $\tau$, in terms of the asymptotic constant, may depend on the smoothness class under consideration. The above result shows that, at least not over a wide range of classes, it does not. Furthermore, Loh and Zhang [17] proposed an adaptive method to determine $\tau$ in the formula $m_n = \tau \log n / \log \log n$ for fixed $n$, but the behavior of this adaptive method was not proved to be better than for $\tau$ constant. Corollary 2 shows that such an adaptive procedure is not needed and that $\tau$ can be taken equal to one.

The final part of Corollary 2, saying that the projection estimator is asymptotically minimax efficient in the Poisson case, is a theoretical argument corroborating the conclusions of the empirical study of Loh and Zhang [17]. Indeed, in a simulation study they compared the projection estimator to a kernel estimator and found that the former performed significantly better for finite sample sizes. Both estimators achieve the optimal rate, but the kernel estimator does not exploit the polynomial structure of the classes $\mathcal{C}(\mathbf{u}^\alpha, C, r)$; this probably introduces a nonnegligible constant in its asymptotics.

5.4. *Proof of Theorem* 3. Throughout the proof we denote by $K_i$ some constants depending only on $(a_k)$, $a$ and $b$.

We start by proving (a). Take $f$ in $\mathcal{C}(\mathbf{u}^\alpha, C, r) \cap \mathbb{H}_1$. In the bias-variance decomposition of Proposition 3, the first term is then bounded by $C^2(m+1)^{-2\alpha}$ for all $m \geq r$. We now bound the second term in the right-hand side of (10). Using (34) and (35) and recalling that the $q_k^\nu$ are orthonormal in $\mathbb{H}$, $\mathrm{var}_f(\widehat{f}_{m,1})$ is bounded by

$$
\mathbb{E}_f \|\check{f}_{m,1}\|_{\mathbb{H}}^2 \leq Z(b)^2 \mathbb{E}_f \left\| \sum_{k=0}^{m-1} \frac{Q_{k,X_1}^\nu}{a_{X_1}} q_k^\nu \right\|_{\mathbb{H}}^2
$$
$$
= Z(b)^2 \sum_{0 \leq l \leq k < m} \left( \frac{Q_{k,l}^\nu}{a_l} \right)^2 \pi_f \mathbf{1}_l.
$$
(50)

Moreover,

$$
\pi_f \mathbf{1}_l = \int_a^b f(\theta) a_l \theta^l \widetilde{Z}(\theta) \, d\theta \leq a_l b^l \widetilde{Z}(a).
$$

Finally, using the bound on $\sum_l (Q_{k,l}^\nu)^2$ given by Lemma A.1, we obtain, for $\lambda_0 > \lambda$ (recall that $\lambda > 1$ depends only on $a$ and $b$) and $m \geq r$,

$$
\mathbb{E}_f \|\widehat{f}_{m,n} - f\|_{\mathbb{H}}^2 \leq C^2 (m+1)^{-2\alpha} + \frac{K_1 \lambda_0^{2m}}{n} \max_{0 \leq k < m} \frac{b^k}{a_k}.
$$

Taking $\lambda_0$ in $(\lambda, \lambda_1^{1/2})$, this proves (a) [see also remark (iii) following the statement of the theorem].

Next we prove (b). Let $f_0 = c\Pi \mathbf{1}_0 = c\pi.(0) = ca_0 \widetilde{Z}$ with $c > 0$ such that $f_0$ is in $\mathbb{H}_1$. Observing that the conditions given in Proposition 6 are satisfied for our choice of $\nu$, we can apply Theorem 1. It remains to verify that (18) implies (47) in this context. Since $f_0$ is in all $V_m$ for positive $m$, $f_0 + \mathcal{C}_{f_0}(K, \mathbf{u}^\alpha, C, r) \subseteq \mathcal{C}(\mathbf{u}^\alpha, C, r)$ for any $K > 0$, say $K = 1$ and any positive $r$. Thus, (18) provides a lower bound on the left-hand side of (47). Next we lower bound the right-hand side of (18). Note that, for this choice of $f_0$, (41) holds for $c_1$ and $c_2$ depending only on $(a_k)$, $a$ and $b$, so that, by (44) and the observation $K_{\infty, f_0}(V_{m+2} \ominus V_m) \leq K_{\infty, f_0}(V_{m+2})$, we find that, for all $m$,

$$(51) \qquad \frac{1}{K_{\infty, f_0}(V_{m+2} \ominus V_m)} \wedge Cu_{m+1} \geq \frac{K_2}{m+2} \wedge C(m+2)^{-\alpha}.$$

Since we have assumed $\alpha > 1$, the right-hand side of this expression equals $C(m+2)^{-\alpha}$ for large $m$. Regarding the second factor in the right-hand side of (18), we note that since $f_0$ is bounded on $[a, b]$, there is a positive constant $K_3$ such that $\pi_{f_0} h \leq K_3 \nu \Pi h$ for all $h \geq 0$. Thus, $\pi_{f_0}\{0, \ldots, m'_n - 1\} = 1 - \sum_{k \geq m'_n} \pi_{f_0} \mathbf{1}_k \geq 1 - K_3 n^{-1} w_n$, so that, under the assumption $w_n \to 0$,

$$\pi^n_{f_0}\{0, \ldots, m'_n - 1\} \geq \exp(n \log(1 - K_3 n^{-1} w_n)) \sim \exp(-K_3 w_n).$$

By applying Theorem 1 as explained above, the two last displayed equations prove (47), with the more detailed lower bound claimed in remark (iii) following the statement of the theorem.

Finally we show (c). The rate of the projection estimator follows from part (a) already proved, and the equivalence relationship (38) between the classes $\mathcal{C}(\mathbf{u}^\alpha, C)$ and $\widetilde{\mathcal{C}}(\alpha, C)$. We now turn to proving optimality of this rate for $\alpha \geq 1$. Optimality over $\widetilde{\mathcal{C}}(\alpha, C)$ is established as in the proof of (b), but taking $f_0(t) = 1/(b-a)$ for $a \leq t \leq b$: we apply Theorem 1 and verify that, for this choice of $f_0$, (18) implies the lower bound

$$\inf_{\hat{f} \in \mathcal{S}_n} \sup_{f \in \widetilde{\mathcal{C}}(\alpha, C) \cap \mathbb{H}_1} \mathbb{E}_f \|\hat{f} - f\|_{\mathbb{H}}^2 \geq K_4 m_n^{-2\alpha}$$

for large enough $n$. Here are the details. Since $f_0$ is in $\widetilde{\mathcal{C}}(\alpha, 1/\sqrt{b-a})$ (as pointed out in Lemma 3) and $\omega_r^\phi(f, t)_p$ is a semi-norm in $f$, $f_0 + \widetilde{\mathcal{C}}(\alpha, \delta) \subseteq \widetilde{\mathcal{C}}(\alpha, C)$ whenever $\delta + 1/\sqrt{b-a} \leq C$. By Proposition 7, $f_0 + \mathcal{C}(\mathbf{u}^\alpha, C_1 \delta) \subseteq f_0 + \widetilde{\mathcal{C}}(\alpha, \delta)$ for all $\delta > 0$, where $C_1$ is as in (38). Adding a constraint on the supremum norm makes an even smaller class, whence $f_0 + \mathcal{C}_{f_0}(K, \mathbf{u}^\alpha, C_1 \delta) \subseteq \widetilde{\mathcal{C}}(\alpha, C)$ for any $K > 0$, say $K = 1$, provided $\delta$ is sufficiently small. Thus, the lower bound of Theorem 1 applies. Using the same arguments as in the proof of (b), we find that this lower bound behaves as $K_5((m'_n + 2)^{-2\alpha} \wedge (m'_n + 2)^{-2})$, which has same rate as $m_n^{-2\alpha}$ for all $\alpha \geq 1$. This completes the proof of (c).

5.5. *Proof of Proposition* 8. We first note that, for $R < \infty$, the proof is completely contained in the proof of Corollary 1. Thus, we here consider the case $R = \infty$ only. Then the conclusions drawn from the subadditive lemma [the paragraph following (A2)] are not helpful, as we can only conclude that $a_k = O(\varepsilon^k)$ for any $\varepsilon > 0$; the latter is indeed implied by $R = \infty$ alone. Thus, a more refined analysis is necessary, as in the proof of Corollary 2.

First we note that, since $b < R$, it holds that $a_k b^k = O(\varepsilon^k)$ for any $\varepsilon > 0$. Thus,

$$\max_{0 \leq k < m} \frac{b^k}{a_k} = \max_{0 \leq k < m} \frac{a_k b^k}{a_k^2} = O\left(\frac{1}{\min_{0 \leq k < m} a_k^2}\right)$$

and, consequently, the condition (45) on $(m_n)$ is implied by rather requiring

$$(52) \qquad \frac{1}{n} \frac{\lambda_1^{2m_n}}{\min_{0 \leq k < m_n} a_k^2} \to 0.$$

The reason why (52) is not used in Theorem 3 is that one would lose the constant derived for the Poisson case in Corollary 2. Moreover, using Lemma 4, we see that for $(w_n)$ to converge to zero with $m'_n = [\eta m_n]$, it is sufficient that

$$(53) \qquad n a_{[\eta m_n]} b^{\eta m_n} \to 0.$$

It remains to check that there exist $\lambda_1 > \lambda$, $\eta > 0$ and $(m_n)$ such that (52) and (53) hold true. We will do this by a constructive proof.

To this end, take $\lambda_1 > \lambda$ arbitrary. The cornerstone in the construction is the following claim, to be proved below: we can find $\eta$, a positive number $C_1$ and $K$ in $(0,1)$ such that, for all $m \geq 0$,

$$(54) \qquad \frac{\lambda_1^{2(m+1)}}{\min_{0 \leq k < m+1} a_k^2} \leq C_1 K^m \frac{b^{-\eta m}}{a_{[\eta m]}}.$$

Given that (54) holds, put

$$m_n = \max\left\{m : \frac{\lambda_1^{2m}}{\min_{0 \leq k < m} a_k^2} K^{-m/2} \leq n\right\}.$$

Since $\lambda_1 > \lambda > 1$ and $K < 1$, the sequence $\lambda_1^{2m} / \min_{0 \leq k < m} a_k^2 \times K^{-m/2}$ is nondecreasing in $m$ and tends to infinity. Thus, $m_n$ is finite for all $n$ and $(m_n)$ is nondecreasing and tends to infinity. Moreover, $\lambda_1^{2m_n} / \min_{0 \leq k < m_n} a_k^2 \leq n K^{m_n/2}$, and (52) follows. On the other hand, from the definition of $m_n$ and (54),

$$n < \frac{\lambda_1^{2(m_n+1)}}{\min_{0 \leq k < m_n+1} a_k^2} K^{-(m_n+1)/2} \leq C_1 K^{(m_n-1)/2} \frac{b^{-\eta m_n}}{a_{[\eta m_n]}}$$

for large $n$. Hence, (53) follows and the construction is complete.

It now remains to prove the claim (54). To do this, let $(r_m)_{m \geq 0}$ be a sequence such that $a_{r_m} = \min_{0 \leq k < m} a_k$. Applying (A2) yields

$$c_0 \min_{0 \leq k < m} a_k^2 = c_0 a_{r_m}^2 \geq a_{2r_m} \qquad \text{for all } m \geq 0.$$

Now put $s = 2r_{m+1}$ and $t = [\eta m]$ and note that since $2r_{m+1} \leq 2m$, for any $\eta \geq 2$, it holds that $t \geq s$. In addition, fix $p$ sufficiently large that $c_0 a_p < 1$ and $c_0 a_p b^p < 1$. The latter can be done since, as noted above, $a_k b^k = O(\varepsilon^k)$ for any $\varepsilon > 0$. Applying (A2) repeatedly, we easily obtain that there is a constant $c_1 > 0$ such that

$$a_t \leq c_1 a_s (c_0 a_p)^{[(t-s)/p]}.$$

Using $2r_{m+1} \leq 2m$ again, we find that $[(t-s)/p] \geq (t-s)/p - 1 \geq \{(\eta - 2)m - 1\}/p - 1$. Recalling that $c_0 a_p < 1$, together with the two last displays, this yields

$$a_{[\eta m]} \leq c_0 c_1 \left( \min_{0 \leq k < m+1} a_k^2 \right) (c_0 a_p)^{((\eta - 2)m - 1)/p - 1}$$

$$\leq C_0 \left( \min_{0 \leq k < m+1} a_k^2 \right) \left( \frac{(c_0 a_p)^{\eta/p}}{(c_0 a_p)^{2/p}} \right)^m$$

for a positive $C_0$. Then with

$$K := \lambda_1^2 \frac{((c_0 a_p)^{1/p} b)^\eta}{(c_0 a_p)^{2/p}},$$

we see that (54) holds for a positive $C_1$. Finally, since $(c_0 a_p)^{1/p} b < 1$, we can choose a large enough $\eta$ such that $K < 1$. The proof is complete.

**6. Discrete deconvolution.** In this section we take $\mathcal{X} = \Theta = \mathbb{Z}$ and let $\zeta$ be counting measure. Let also $p$ be a fixed and known probability mass function on $\mathbb{Z}$ and consider the mixands $\pi_\theta(\cdot) = p(\cdot - \theta)$. Another way to view this setup is the following. Take independent random variables $\varepsilon$ and $U$, both in $\mathbb{Z}$, with probability mass functions $p$ and $f$, respectively, and put $X = U + \varepsilon$. Then the probability mass function of $X$ is the convolution $(f \star p)(\cdot) = \sum_\theta p(\cdot - \theta) f(\theta)$ of $f$ and $p$, which we can also write as $\sum_\theta \pi_\theta(\cdot) f(\theta)$. Our interest in recovering $f$ from i.i.d. observations from $X$ can thus be phrased as a deconvolution problem. Note that this setting includes the case of $\varepsilon$ being zero, that is, we estimate a discrete distribution.

Observe that, for all integer $k$, $\Pi \mathbf{1}_k = p(k - \cdot)$. Applying the general approach of Section 4, we take $V_0 = \{0\}$ and, for all $m \geq 1$, $V_m := \text{Span}(\Pi \mathbf{1}_k : |k| < m)$. This, of course, defines an increasing sequence of linear spaces. It remains to choose the measure $\nu$ or, equivalently, the space $\mathbb{H} = L^2(\nu)$. A natural choice is to let $\nu$ be counting measure, that is, $\mathbb{H} = l^2(\mathbb{Z})$. Then, since $p$ is square-summable, $(V_m)$ is a sequence of subspaces of $\mathbb{H}$. It is practical to define the projection estimator using Fourier series. Thus, let

$$p^*(\lambda) = \sum_{k \in \mathbb{Z}} p(k) e^{-ik\lambda} \qquad \text{for all } \lambda \in (-\pi, \pi]$$

be the Fourier series with coefficients $(p(k))_{k \in \mathbb{Z}}$. Then $p^* \in L^2(-\pi, \pi]$ and, because $p$ is a density, $p^*$ is continuous with positive $L^2$ norm. The Fourier series with coefficients $(\Pi \mathbf{1}_k)_{k \in \mathbb{Z}}$ simply reads $p^*(-\lambda) e^{-ik\lambda}$. Because there necessarily is an interval on which $p^*$ is nonzero, $\{p^*(-\lambda) e^{-ik\lambda}\}_{k \in \mathbb{Z}}$ is linearly independent and assumption (A1) then follows immediately. Hence, the projection estimator is well defined and the results of Section 4 apply.

Let us derive the expression for the projection estimator $\widehat{f}_{m,n}$ through the Fourier series $\widehat{f}^*_{m,n}$ with Fourier coefficients $(\widehat{f}_{m,n}(k))_{k \in \mathbb{Z}}$. Let $P^*_n$ be the Fourier series associated to the coefficients $(P_n \mathbf{1}_k)_{k \in \mathbb{Z}}$,

$$P^*_n(\lambda) = \sum_{k \in \mathbb{Z}} (P_n \mathbf{1}_k) e^{-ik\lambda} \qquad \text{for all } \lambda \in (-\pi, \pi].$$

Then applying Parseval's formula to (7) with $g = \Pi \mathbf{1}_k$, $\widehat{f}^*_{m,n}$ is the unique element in $\mathrm{Span}(e^{-ik\lambda} p^*(-\lambda) : |k| < m)$ which satisfies, for all $k = -m, \ldots, m$,

$$(55) \qquad \frac{1}{2\pi} \int_{-\pi}^{\pi} \widehat{f}^*_{m,n}(\lambda) p^*(\lambda) e^{ik\lambda} \, d\lambda = P_n \mathbf{1}_k = \frac{1}{2\pi} \int_{-\pi}^{\pi} P^*_n(\lambda) e^{ik\lambda} \, d\lambda.$$

Here we will treat the special case where the following condition holds:

$$(56) \qquad K_p := \int_{-\pi}^{\pi} \frac{1}{|p^*(\lambda)|^2} \, d\lambda < \infty.$$

This condition implies that $p^*$ may only vanish on a Lebesgue null set, and in a singular way (it cannot have a finite derivative where it vanishes). It is, of course, verified, for instance, if $|p^*|$ is bounded away from zero, which includes the case of estimating a discrete distribution ($\varepsilon = 0$). If (56) holds, there is a function in $L^2(-\pi, \pi]$ such that (55) holds for all $k \in \mathbb{Z}$. Indeed, since $P_n$ is a probability, $|P^*_n|$ is bounded by one so that $P^*_n / p^*$ is in $L^2(-\pi, \pi]$ and satisfies (55) for all $k$ whenever (56) holds. This amounts to taking the definition of the projection estimator to its limit $m = \infty$; hence, we put $\widehat{f}^*_{\infty,n} := P^*_n / p^*$ or, equivalently,

$$\widehat{f}_{\infty,n}(k) := \frac{1}{2\pi} \int_{-\pi}^{\pi} \frac{P^*_n(\lambda)}{p^*(\lambda)} e^{ik\lambda} \, d\lambda \qquad \text{for all } k \in \mathbb{Z}.$$

To compute the expectation of $\widehat{f}_{\infty,n}$, first apply (8) and then Parseval's formula to obtain

$$\mathbb{E}_f P^*_n(\lambda) = \sum_{k \in \mathbb{Z}} \pi_f \mathbf{1}_k e^{-ik\lambda}$$

$$= \sum_{k \in \mathbb{Z}} (f, \Pi \mathbf{1}_k)_{\mathbb{H}} e^{-ik\lambda} = f^*(\lambda) p^*(\lambda) \qquad \text{for all } \lambda \text{ in } (-\pi, \pi],$$

where $f^*$ is the Fourier series with coefficients $(f(k))_{k \in \mathbb{Z}}$. The two last equations give

$$(57) \quad \mathbb{E}_f \widehat{f}_{\infty,n}(k) = \frac{1}{2\pi} \int_{-\pi}^{\pi} \frac{\mathbb{E}_f P_n^*(\lambda)}{p^*(\lambda)} e^{ik\lambda} \, d\lambda = f(k) \qquad \text{for all } k \in \mathbb{Z},$$

where, as $|P_n^*/p^*| \leq 1/|p^*|$ is integrable, we applied Fubini's theorem. Hence, $\widehat{f}_{\infty,n}$ is unbiased, as expected, since it corresponds to $m = \infty$. Furthermore,

$$(58) \quad \operatorname{var}_f(\widehat{f}_{\infty,n}) = \frac{1}{n} \operatorname{var}_f(\widehat{f}_{\infty,1}) \leq \frac{1}{2\pi n} \mathbb{E}_f \int_{-\pi}^{\pi} |\widehat{f}_{\infty,1}^*(\lambda)|^2 \, d\lambda \leq \frac{K_p}{2\pi n}.$$

As shown by the following result, we are in a case where the minimax MISE rate is achieved by $\widehat{f}_{\infty,n}$ over any reasonable subclass of $\mathbb{H}_1$. This is a degenerate case compared to the general setting of Section 4, because smoothness conditions of the form (14) do not improve the minimax rate. To formulate the result, we define the line segment $[f_0, f_1]$ between $f_0$ and $f_1$ as $[f_0, f_1] := \{(1-w)f_0 + wf_1 : w \in [0, 1]\}$.

THEOREM 4. *Assume that* (56) *holds. Then the MISE of the projection estimator $\widehat{f}_{\infty,n}$ has rate $n^{-1}$ over $\mathbb{H}_1$, and this rate is minimax over any class $\mathcal{C}$ included in $\mathbb{H}_1$ and containing a line segment $[f_0, f_1]$ for distinct $f_0$ and $f_1$ in $\mathbb{H}_1$. More precisely, for any $f_0$ and $f_1$ in $\mathbb{H}_1$ and any positive integer $n$,*

$$c_0 (K_1 + K_2 c_1 n)^{-1} \leq \inf_{\hat{f} \in \mathcal{S}_n} \sup_{f \in [f_0, f_1]} \mathbb{E}_f \|\hat{f} - f\|_{\mathbb{H}}^2$$

$$\leq \sup_{f \in \mathbb{H}_1} \mathbb{E}_f \|\widehat{f}_{\infty,n} - f\|_{\mathbb{H}}^2 \leq \frac{K_p}{2\pi n},$$

*where $K_1$ and $K_2$ are universal positive constants, $c_0 := \sum_{l \in \mathbb{Z}} (f_1(l) - f_0(l))^2$ and $c_1 := \sum_{l \in \mathbb{Z}} |\pi_{f_1}(l) - \pi_{f_0}(l)|$.*

REMARK. Observe that the lower bound may reduce to a positive constant only when the model is not identifiable, that is, if there exist $f_0$ and $f_1$ in $\mathbb{H}_1$ such that $c_1 = 0$. In this case (56) cannot be fulfilled.

PROOF OF THEOREM 4. Since $\widehat{f}_{\infty,n}$ is unbiased, the upper bound simply is (58).

The rate $n^{-1}$ of the lower bound on the MISE generally holds for any regular parametric statistical model. Here we derive it via a Bayes risk lower bound. Consider the parametric model $\{\pi_{(1-w)f_0 + wf_1} : w \in (0, 1)\}$. Put any continuously differentiable prior density $r$ on $w \in (0,1)$ which is symmetric about $1/2$: $r(1/2 + w) = r(1/2 - w)$. Let $\mathcal{I}(w) = \mathbb{E}_{(1-w)f_0+wf_1} \times [(\partial_w \log \pi_{(1-w)f_0+wf_1}(X))^2]$ and $\mathcal{I}(r) = \int_0^1 \dot{r}(w)^2/r(w)\,dw$. Then, for any estimator $\hat{f}$ from $n$ observations,

$$\sup_{f \in [f_0, f_1]} \mathbb{E}_f \|\hat{f} - f\|_{\mathbb{H}}^2$$

(59)
$$\geq \int_0^1 \mathbb{E}_{(1-w)f_0+wf_1} \|\hat{f} - ((1-w)f_0 + wf_1)\|_{\mathbb{H}}^2 r(w)\,dw$$

$$= \int_0^1 \sum_{k \in \mathbb{Z}} \mathbb{E}_{(1-w)f_0+wf_1}[(\hat{f}(k) - ((1-w)f_0(k) + wf_1(k)))^2] r(w)\,dw$$

$$\geq \sum_{k \in \mathbb{Z}} \frac{(f_1(k) - f_0(k))^2}{n \int_0^1 \mathcal{I}(w)r(w)\,dw + \mathcal{I}(r)}.$$

Here the first inequality is the Bayes risk lower bound on the minimax risk and the second one is the van Trees inequality (a Bayesian Cramér–Rao bound); see, for example, [11]. We easily compute

$$\mathcal{I}(w) = \sum_{k \in \mathbb{Z}} \frac{(\pi_{f_1}(k) - \pi_{f_0}(k))^2}{(\pi_{f_0}(k) + \pi_{f_1}(k))/2 + (w - 1/2)(\pi_{f_1}(k) - \pi_{f_0}(k))},$$

with the convention $0/0 = 0$. Using the symmetry of $r$, we obtain

$$\int_0^1 \mathcal{I}(w)r(w)\,dw$$

$$= \sum_{k \in \mathbb{Z}} \int_0^1 \frac{|\pi_{f_1}(k) - \pi_{f_0}(k)|r(w)\,dw}{(\pi_{f_0}(k) + \pi_{f_1}(k))/(2|\pi_{f_1}(k) - \pi_{f_0}(k)|) + (w - 1/2)}$$

$$\leq \sum_{k \in \mathbb{Z}} |\pi_{f_1}(k) - \pi_{f_0}(k)| \int_0^1 w^{-1} r(w)\,dw,$$

where we simply used $a + b \geq |a - b|$ for any nonnegative $a$ and $b$ for the inequality. Hence, (59) gives the required lower bound, where $K_1 = \mathcal{I}(r)$ and $K_2 = \int_0^1 w^{-1} r(w)\,dw$ are fixed once a particular choice of $r$ is made. □

REMARK. There is a tradeoff between $K_1$ and $K_2$ when the prior density $r$ is chosen for optimizing the lower bound for finite $n$. Asymptotically, the lower bound is equivalent to $c_0/(K_2 c_1 n)$, and it is easy to see that the infimum of possible values of $K_2$ is 2 (let $r$ tend to a point mass located at $w = 1/2$). Hence,

$$\liminf_{n \to \infty} \inf_{\hat{f} \in \mathcal{S}_n} \sup_{f \in [f_0, f_1]} n \mathbb{E}_f \|\hat{f} - f\|_{\mathbb{H}}^2 \geq \frac{c_0}{2c_1}.$$

The right-hand side depends on $f_0$ and $f_1$ and should be compared to the asymptotic upper bound $K_p/(2\pi)$.

The statistical literature on deconvolution is vast, but it is primarily concerned with continuous random variables having densities with respect to Lebesgue measure, often on $\mathbb{R}$. Some key references on achievable minimax rates of convergence over suitable smoothness classes in that setting are [2] and [6] for pointwise estimation of the mixing density, and [5, 7, 8, 22] for (weighted) $L^p$ loss. The difficulty of the estimation then depends on whether the characteristic function of $\varepsilon$, that is, essentially our $p^*$, vanishes algebraically or exponentially fast at infinity, these cases being referred to as ordinary smooth and supersmooth error densities, respectively. With an ordinary smooth error density, the optimal rate of convergence is algebraic in $n$, whereas it is algebraic in $\log n$ when the error density is supersmooth.

In the discrete setting considered here, the notions of ordinary smooth and supersmooth error densities are void, since $p^*$ is defined on a compact interval, the unit circle. The MISE rate $n^{-1}$ of Theorem 4 is also faster than what is obtained in the papers cited above; it only appears as a limit in the ordinary smooth case when the unknown density $f$ has infinite smoothness. In the discrete case, the rate $n^{-1}$ may not hold when (56) fails; some additional remarks on this issue are given in Section 8.

**7. Mixtures of uniform discrete distributions.** We now take $\Theta = \mathbb{N} := \{1, 2, \ldots\}$, $\mathcal{X} = \mathbb{Z}_+$ and let $\zeta$ and $\nu$ both be counting measure. Thus, $\mathbb{H} = l^2(\Theta)$. Consider mixands given by the family of uniform discrete distributions on $\{0, 1, \ldots, \theta - 1\}$; that is, $\Pi \mathbf{1}_k(\theta) = \pi_\theta(k) = \theta^{-1}$ for $0 \leq k \leq \theta - 1$ and 0 otherwise. Observe that, for all $k \geq 0$, $\Pi \mathbf{1}_k(\theta) = \theta^{-1}$ for large $\theta$ so that $\Pi \mathbf{1}_k$ is not in $l^1(\Theta)$ and (A1) does not hold. Then letting the space $V_m$ be spanned by $(\Pi \mathbf{1}_k)_{0 \leq k < m}$ as in Section 4 would yield an estimator $\hat{f}_{m,n}$ that is a linear combination of nonintegrable functions and, hence, a poor estimator of the mixing density. It is possible to circumvent this problem by replacing $\nu$ by a distribution such that $((1+\theta)^{-1})_{\theta \geq 0}$ belongs to $L^1(\nu)$, but then the difficulty would lie in the definition of the smoothness classes (14). Indeed, in this case a different choice of $V_m$ provides a much simpler definition of smoothness classes. For all $k \geq 0$, we let $h_k = (k+1)(\mathbf{1}_k - \mathbf{1}_{k+1})$, which yields

$$\Pi h_k = (k+1)(\Pi \mathbf{1}_k - \Pi \mathbf{1}_{k+1}) = \mathbf{1}_{k+1}.$$

Hence, $(\Pi h_k)_{k \geq 0}$ is itself the orthonormal basis denoted by $(\phi_k)$ in Section 3. It follows that

$$V_m = \mathrm{Span}(\Pi h_k, 0 \leq k < m) = \{f \in l^2(\Theta) : f(\theta) = 0 \text{ for all } \theta > m\},$$

the projection estimator is

$$\widehat{f}_{m,n}(k) = k P_n(\mathbf{1}_{k-1} - \mathbf{1}_k)\mathbf{1}(k \leq m) \qquad \text{for all } k \geq 1,$$

and the smoothness classes of Section 4 read

$$
(60) \qquad \mathcal{C}(u, C, r) = \left\{ f \in l^2(\Theta) \colon \sum_{k > m} f^2(k) \leq C^2 u_m^2 \text{ for all } m \geq r \right\}.
$$

Since $\operatorname{var}_f(h_k) \leq (k+1)^2(\pi_f(k) + \pi_f(k+1))$ for any $f$ in $\mathbb{H}_1$, Proposition 3, along with (11), gives, for all $m$,

$$
(61) \quad \mathbb{E}_f \|\hat{f}_{m,n} - f\|_\mathbb{H}^2 \leq \sum_{\theta > m} f^2(\theta) + \frac{1}{n} \sum_{k=0}^{m-1} (k+1)^2 (\pi_f(k) + \pi_f(k+1)).
$$

As in Section 6, this implies that the MISE rate $n^{-1}$ is achievable as soon as $\pi_f$ has finite second moment, that is, $\sum_{k \geq 0} k^2 \pi_f(k) < \infty$. Moreover, it holds that $\pi_f$ has finite second moment for any $f$ in $\mathcal{C}(u, C, r) \cap \mathbb{H}_1$ whenever $u = (u_m)$ satisfies $\sum m^{3/2} u_m < \infty$. Indeed, as a simple consequence of the Cauchy–Schwarz inequality, $\pi_f(k) = O(\{k^{-1} \sum_{\theta > k} f^2(\theta)\}^{1/2}) = O(k^{-1/2} u_k)$ for such $f$. The interesting cases are thus those when $(u_m)$ decreases slowly, and Corollary 3 below provides such an example.

THEOREM 5. *Let $u = (u_m)$ be a positive sequence decreasing to zero, $C$ a positive number and $r$ a positive integer. Then for all sufficiently large $m$,*

$$
(62) \quad \inf_{\hat{f} \in \mathcal{S}_n} \sup_{f \in \mathcal{C}(u, C, r) \cap \mathbb{H}_1} \mathbb{E}_f \|\hat{f} - f\|_\mathbb{H}^2 \geq \left( \frac{C u_{m+1}}{2} \right)^2 \left( 1 - \frac{\sqrt{5}}{2m} C u_{m+1} \right)^n,
$$

*and for any integer $m \geq r$,*

$$
(63) \qquad \sup_{f \in \mathcal{C}(u, C, r) \cap \mathbb{H}_1} \mathbb{E}_f \|\hat{f}_{m,n} - f\|_\mathbb{H}^2 \leq (C u_m)^2 + \frac{2m^2}{n}.
$$

PROOF. The upper bound (63) follows from (61) and the bound $\sum_{k=0}^{m-1} (k+1)^2 \pi_f(k) \leq m^2$.

The lower bound (62) is obtained along the same lines as the lower bound of Theorem 1. We apply Proposition 2 with $h(k) = \mathbf{1}(k < m)$,

$V = \operatorname{Span}\{h(k)\pi.(k), k \in \mathcal{X}\}$

$= \{f \in l^2(\Theta) \colon \text{there exists } \lambda \in \mathbb{R} \text{ such that } f(\theta) = \lambda \theta^{-1} \text{ for all } \theta \geq m\}$

and $\mathcal{C}^\star := V^\perp \cap \mathcal{C}(u, \frac{1}{2}C, r)$. Thus, the assumptions of Proposition 2 are immediately satisfied. To arrive at (62), we still need to find a suitable probability density $f_0$ such that $f_0 + \mathcal{C}^\star \subseteq \mathcal{C}(u, C, r)$, and a function $g$ in $\mathcal{C}^\star$ such that $f_0 \pm g$ are in $\mathbb{H}_1$ to bound the supremum in the lower bound of Proposition 2 from below.

In order to do this, let $g$ in $l^2(\Theta)$ be one of the two sequences satisfying the three equations

(i) $g(\theta) = 0$ unless $\theta = m$, $m+1$ or $m+2$;

(ii) $\sum_{\theta=m}^{m+2} g(\theta) = 0$ and $\sum_{\theta=m}^{m+2} \theta^{-1} g(\theta) = 0$;

(iii) $\sum_{\theta=m}^{m+2} g^2(\theta) = (\frac{1}{2} C u_{m+1})^2$.

Then $\|g\|_{\mathbb{H}} = \frac{1}{2} C u_{m+1}$.

Let $f_0$ be such that $f_0(\theta) = |g(\theta)|$ for $\theta > 1$ and $f_0(1) = 1 - \sum_{\theta=m}^{m+2} f_0(\theta)$. Then for $m$ sufficiently large, $\sum_{\theta=m}^{m+2} f_0(\theta) \leq \sqrt{3} (\sum_{\theta=m}^{m+2} g^2(\theta))^{1/2} = \frac{\sqrt{3}}{2} C u_{m+1}$ is less than one. Hence, $f_0$ belongs to $\mathbb{H}_1$ for large $m$. Using (60) and the result that $(u_k)$ is nonincreasing, one readily checks that $f_0$ is in $\mathcal{C}(u, \frac{1}{2}C, r)$. It is then immediate that $f_0 + \mathcal{C}(u, \frac{1}{2}C, r) \subseteq \mathcal{C}(u, C, r)$, and thus also $f_0 + \mathcal{C}^\star \subseteq \mathcal{C}(u, C, r)$.

We now proceed to checking that $g$ belongs to $\mathcal{C}^\star$ and that $f_0 \pm g$ are in $\mathbb{H}_1$. The latter follows from $\sum g(\theta) = 0$ and $|f_0| \geq |g|$. The former is also true as $g$ both belongs to $\mathcal{C}(u, \frac{1}{2}C, r)$, which is checked as for $f_0$, and is perpendicular to $V$, which follows from item (i) and the second part of item (ii) in its definition.

Hence, Proposition 2 gives

$$\inf_{\hat{f} \in \mathcal{S}_n} \sup_{f \in \mathcal{C}(u,C,r) \cap \mathbb{H}_1} \mathbb{E}_f \|\hat{f} - f\|_{\mathbb{H}}^2 \geq \|g\|_{\mathbb{H}}^2 \pi_{f_0}^n \{0, 1, \ldots, m-1\}.$$

We easily compute

$$\pi_{f_0}\{0, 1, \ldots, m-1\}$$
$$= 1 - f_0(m+1)(m+1)^{-1} - 2 f_0(m+2)(m+2)^{-1}$$
$$\geq 1 - \frac{\sqrt{5}}{m} \|g\|_{\mathbb{H}},$$

and (62) follows from the two last inequalities. □

COROLLARY 3. *Let $\alpha$ and $C$ be two positive numbers and $r$ a nonnegative integer. Then the minimax MISE rate over the class $\mathcal{C}((\log^{-\alpha}(1+n)), C, r)$ is $(\log n)^{-2\alpha}$. This rate is achieved by the projection estimator $\hat{f}_{m_n, n}$ with $m_n = [\tau n^\beta]$ for any positive number $\tau$ and any positive $\beta$ less than $1/2$. Moreover, this estimator is asymptotically MISE efficient up to a factor $4/\beta^{2\alpha}$.*

PROOF. Put $u_m = (\log(1+m))^{-\alpha}$. Use the lower bound (62) with $m = n$ to obtain the asymptotics $(Cu_n/2)^2$. Then use the upper bound (63), with $m_n = [\tau n^\beta]$ and $\beta \in (0, 1/2)$, to obtain

$$u_{m_n}^{-2} \sup_{f \in \mathcal{C}(u,C,r) \cap \mathbb{H}_1} \mathbb{E}_f \|\hat{f}_{m_n, n} - f\|_{\mathbb{H}}^2 \leq C^2 + \frac{2 m_n^2}{n u_{m_n}^2} = C^2 + o(1).$$

Finally, $(Cu_{m_n})^2/(Cu_n^2/2)^2 \to 4/\beta^{2\alpha}$. The proof is complete. □

**8. Open problems.** In Sections 5 and 6 we have investigated some particular cases for which (A1) is satisfied and, thus, both Theorem 1 and Theorem 2 apply. The upper bounds were, in fact, obtained directly without using Theorem 2, but are essentially similar (see the remark following the proof of Theorem 2). We have also seen in Section 7 that, in certain situations when (A1) does not hold, this approach could be adapted. However, only in the case of power series mixtures and in the adaptation of Section 7 could we compute explicit lower and upper bounds giving rise to identical rates. In the case of discrete deconvolution of Section 6, we even gave an alternative lower bound which gives the optimal rate in the degenerate case where the projection estimator $\widehat{f}_{m,n}$ can be defined for $m = \infty$. In this section we outline a few open problems for which the framework of the present paper is potentially applicable, but in which we do not attempt to compute the bounds.

8.1. *Multivariate power series.* Proposition 6 shows that Theorem 1 applies for a large range of dominating measures $\nu$. Condition (i) is a simple reformulation of the requirement of having $\pi.(k)$ both in $L^1(\nu)$ and in $L^2(\nu)$ for all $k$, and easily generalizes to other mixands. Condition (ii) also generalizes when $\pi.(k)$ is related to a well-known sequence of linearly independent functions (here the polynomials). For instance, it trivially generalizes to multivariate power series mixands (see [13], Chapter 38). A slightly different setting concerns the bivariate Poisson distributions ([13], Chapter 37) given by, for all $\theta = (\theta_1, \theta_2, \theta_{12})$ in $\Theta := (0, \infty)^3$ and all nonnegative integers $x_1$ and $x_2$,

$$\pi_\theta(x_1, x_2) = e^{-(\theta_1+\theta_2+\theta_{12})} \sum_{i=0}^{x_1 \wedge x_2} \frac{\theta_1^{x_1-i}\theta_2^{x_2-i}\theta_{12}^i}{(x_1-i)!(x_2-i)!i!}.$$

In this case (i) is modified to $\int_\Theta |\theta|^k e^{-(\theta_1+\theta_2+\theta_{12})}\nu(d\theta) < \infty$ and (ii) is unchanged. This is easily seen upon observing that $\{e^{\theta_1+\theta_2+\theta_{12}}(\pi_\theta(x_1,x_2)\}_{x_1,x_2\geq 0}$ is a collection of trivariate polynomials in $\theta = (\theta_1, \theta_2, \theta_{12})$ that are linearly independent as the term $\theta_1^{x_1}\theta_2^{x_2}$ only appears in $\pi_\theta(x_1, x_2)$. The rest of the proof above applies similarly.

Although Theorems 1 and 2 apply, the rates they provide are not known explicitly.

8.2. *Power series mixing distributions with noncompact support.* Let again $\pi_\theta$ be the Poisson distribution with mean $\theta$, and let $\Theta = \mathbb{R}_+$. Thus, we have power series mixands with $a_k = 1/k!$, but with $\Theta$ being unbounded. Take $\nu$ as Lebesgue measure, so that $\mathbb{H} = L^2(\mathbb{R}_+)$. Applying (12) [see also (35)], we find that, for all $f$ in $\mathbb{H}_1$, the variance term in Proposition 3 is bounded by

$$(64) \qquad \mathbb{E}_f \|\widehat{f}_{m,1}\|_{\mathbb{H}}^2 = \mathbb{E}_f \left( \sum_{k=0}^{m-1} \Phi_{k,X_1}^2 \right) \leq \sum_{0 \leq l \leq k < m} \Phi_{k,l}^2,$$

which should be compared to (50). For $\lambda_1 > 2 + \sqrt{17}/2$, this bound, along with Lemma A.2 in the bias-variance decomposition, shows that for any positive number $C$, any nonnegative integer $r$ and any sequences $(u_n)$ and $(m_n)$ satisfying $\lambda_1^{m_n}/u_{m_n} = o(n^{1/2})$,

$$\limsup_{n \to \infty} u_{m_n}^{-2} \sup_{f \in \mathcal{C}(u,C,r) \cap \mathbb{H}_1} \mathbb{E}_f \|\widehat{f}_{u_{m_n},n} - f\|_{\mathbb{H}}^2$$
$$\leq \limsup_{n \to \infty} u_{m_n}^{-2} (C^2 u_{m_n}^2 + n^{-1} K_1 \lambda_1^{2m_n}) = C^2.$$

For sequences $\mathbf{u}^\alpha$ the obtained root MISE rate is $(\log n)^{-\alpha}$ by choosing $m_n = [\tau \log n]$ for small $\tau$. This is better than when $b < \infty$ (see Corollary 2).

Concerning lower bounds on the MISE, Loh and Zhang [17] give such a one in their Theorem 4 over particular classes related to ours, but their assumptions do not apply in the case considered here because they correspond to a weight function $w = \mathbf{1}_{\mathbb{R}_+}$ with infinite $L^1$ norm. Hence, it is still to be found if the logarithmic rate of the projection estimator is optimal in this case. Theorem 1 applies and Proposition 5 shows that $\mathcal{C}(u, C, r) \cap \mathbb{H}_1$ is nonempty for a positive $r$ or for large $C$. The next problem, which we have not solved, rather consists in finding $f_0$ in this intersection such that $u_m = O(1/K_{\infty,f_0}(V_{m+2} \ominus V_m))$, which is the key in our method for showing that the minimax MISE rate is $(u_{m_n}^2)$.

8.3. *Discrete deconvolution with vanishing characteristic function.* Let us return to the setting of Section 6. If condition (56) is not satisfied, our analysis must be refined. It is indeed possible that the optimal rate is slower in this case, and the behavior of $p^*$ at its zeros may then yield the optimal rate of convergence.

To our knowledge, this problem has not been studied. A possible approach would be to mimic that of Section 5.1 by observing that the projection estimator can be easily defined using a sequence of orthonormal trigonometric polynomials in $L^2(\nu')$ with $\nu'(dt) = |p^*|^2(t)\mathbf{1}_{(-\pi,\pi]}(t)\,dt$, and to express $\widehat{f}_{m,n}$ using such a sequence. However, in contrast to power series mixands, the behavior of the projection estimator here should be driven by the measure $\nu'$, that is, by using precise assumptions on its zeros when (56) fails.

## APPENDIX

**Recurrence relations for orthonormal polynomials.** In this appendix we give some further results for the orthogonal polynomials $(q_k^\nu)$ in $\mathbb{H}$ and $(q_k^{\nu'})$ in $\mathbb{H}'$, introduced in Section 5.1. For any measure $\nu_0$ on $\mathbb{R}$, one can construct an orthogonal sequence of polynomials $(r_k)_{k\geq 0}$ with increasing degrees in $L^2(\nu_0)$ using a so-called three terms recurrence relation,

$$r_{k+1}(t) = (t - \alpha_k)r_k(t) - \beta_k r_{k-1}(t) \qquad \text{with } r_{-1} = 0 \text{ and } r_0 = 1,$$

where $(\alpha_k)_{k\geq 0}$ and $(\beta_k)_{k\geq 0}$ are sequences depending on $\nu_0$. Moreover, putting $\beta_0 = \|r_0\|_{\mathbb{H}}^2$, one has ([10], equations (1.13))

$$(65) \qquad N_k := \|r_k\|_{\mathbb{H}} = \left(\prod_{j=0}^{k} \beta_j\right)^{1/2} \qquad \text{for all } k \geq 0.$$

Let the polynomials have coefficients $r_k(t) = \sum_l R_{k,l} t^l$ and put

$$(66) \qquad Q_{k,l}^{\nu_0} = R_{k,l}/N_k \qquad \text{for all } k,l \geq 0.$$

The latter coefficients are those of an orthonormal sequence corresponding to $Q_{k,l}^{\nu'}$ and $Q_{k,l}^{\nu}$ for $\nu_0$ equal to $\nu'$ and $\nu$, respectively. The three terms recurrence relation can be written

$$(67) \quad R_{k+1,l} = R_{k,l-1} - \alpha_k R_{k,l} - \beta_k R_{k-1,l} \qquad \text{for all } k,l \geq 0, \text{ with } R_{0,0} = 1$$

and the convention $R_{k,l} = 0$ if $l < 0$ or $l > k$. Hence, by (65), (67) and (66), we see that knowledge of $(\alpha_k)_{k\geq 0}$ and $(\beta_k)_{k\geq 0}$ provides a simple algorithm for computing the coefficients $Q_{k,l}^{\nu_0}$ recursively at a low computational cost.

Let us now derive the coefficients $\alpha_k$ and $\beta_k$ for particular choices for $\nu_0$.

**A.1. Legendre polynomials.** Let $\nu_0$ be Lebesgue measure on $[-1,1]$. In this case $\alpha_k = 0$, $\beta_0 = 2$ and $\beta_k = 1/(4 - 1/k^2)$ for $k \geq 1$ ([10], equation (2.1)).

**A.2. Translated-scaled Legendre polynomials.** Let $\nu_0$ be Lebesgue measure on an interval $[a,b]$. Denote $\mu := (a+b)/2$ and $\delta := (b-a)/2$. Replacing $t$ by $(u-\mu)/\delta$ in the three terms recurrence relation for Legendre polynomials, one obtains an orthogonal sequence for $\nu_0$ satisfying

$$r_{k+1}\left(\frac{u-\mu}{\delta}\right) = \left(\frac{u-\mu}{\delta} - \alpha_k\right) r_k\left(\frac{u-\mu}{\delta}\right) - \beta_k r_{k-1}\left(\frac{u-\mu}{\delta}\right).$$

Multiplying by $\delta^{k+1}$ and identifying $(\delta^k r_k((u-\mu)/\delta))_k$ with a new orthogonal sequence, the previous equation gives the following coefficients in this case: $\alpha_k = -\mu$, $\beta_0 = 2\delta$ and $\beta_k = \delta^2/(4 - 1/k^2)$ for $k \geq 1$.

The following result serves for bounding the variance of $\check{f}_{m,n}$ [see (34)] when $\nu = \nu_0$.

LEMMA A.1. *Let $\lambda_0$ be a number larger than $\lambda := \gamma + \sqrt{\gamma^2 + 1}$ with $\gamma = (2 + a + b)/(b - a)$. Then $\sum_l (Q_{k,l}^{\nu_0})^2 = O(\lambda_0^{2k})$.*

PROOF. It follows from (67) that
$$\|R_{k+1}\| \le (1 + |\alpha_k|)\|R_k\| + \beta_k \|R_{k-1}\|,$$
where $\|R_k\| := (\sum_l R_{k,l}^2)^{1/2}$. Consequently, dividing by $N_{k+1}$ as given above, we obtain

(68) $$\|Q_{k+1}^{\nu_0}\| \le \frac{1 + |\alpha_k|}{\sqrt{\beta_{k+1}}} \|Q_k^{\nu_0}\| + \sqrt{\frac{\beta_k}{\beta_{k+1}}} \|Q_{k-1}^{\nu_0}\|.$$

Note that
$$\lim_{k \to \infty} \frac{1 + |\alpha_k|}{\sqrt{\beta_{k+1}}} = \frac{2(1 + \mu)}{\delta} \quad \text{and} \quad \sqrt{\frac{\beta_k}{\beta_{k+1}}} \le 1 \quad \text{for all } k \ge 1,$$

and that the positive solution of the quadratic equation

(69) $$x^2 - \frac{2(1+\mu)}{\delta} x - 1 = 0 \quad \text{with } \mu = \frac{a+b}{2} \text{ and } \delta = \frac{b-a}{2}$$

is $\lambda$. The lemma follows. □

**A.3. Laguerre measure.** Let $\nu_0$ be the measure with density $e^{-t}$ with respect to Lebesgue measure over $t \in (0, \infty)$.

In this case $\alpha_k = 2k + 1$, $\beta_0 = 1$ and $\beta_k = k^2$ for $k \ge 1$ ([10], equation (2.4)).

**A.4. Squared Laguerre measure.** Let finally $\nu_0$ be the measure with density $e^{-2t}$ with respect to Lebesgue measure over $t \in (0, \infty)$.

This corresponds to $\nu'$ in the case of Poisson mixands, $a_k = 1/k!$ and $\widetilde{Z}^2(t) = e^{-2t}$. We have a result corresponding to Lemma A.1, but we now study $\mathbb{H}'$ rather than $\mathbb{H}$ as in the case of compact support and Lebesgue measure, as our interest lies in the coefficients $\Phi_{k,l} = Q_{k,l}^{\nu_0}/a_l$ (see Section 5.1).

LEMMA A.2. *Let $\lambda_1$ be a number larger than $2 + \sqrt{17}/2$. Then $\sum_l \Phi_{k,l}^2 = \sum_l (Q_{k,l}^{\nu_0}/a_l)^2 = O(\lambda_1^{2k})$.*

PROOF. Dividing (67) by $a_l$ one obtains

(70) $$\left( \sum_l \frac{Q_{k+1,l}^{\nu_0 2}}{a_l} \right)^{1/2}$$
$$\le \frac{r_k + |\alpha_k|}{\sqrt{\beta_{k+1}}} \left( \sum_l \frac{Q_{k,l}^{\nu_0 2}}{a_l} \right)^{1/2} + \sqrt{\frac{\beta_k}{\beta_{k+1}}} \left( \sum_l \frac{Q_{k-1,l}^{\nu_0 2}}{a_l} \right)^{1/2},$$

where $r_k := \max_{0 \leq l \leq k}(a_{l-1}/a_l)$. By arguments similar to the ones used for the translated-scaled Legendre polynomials, $\alpha_k = k + 1/2$, $\beta_0 = 1/2$ and $\beta_k = k^2/4$ for $k \geq 1$. Since $a_k = 1/k!$, $r_k = k$ in (70). Moreover,

$$\lim_{k \to \infty} \frac{k + |\alpha_k|}{\sqrt{\beta_{k+1}}} = 4 \quad \text{and} \quad \sqrt{\frac{\beta_k}{\beta_{k+1}}} \leq 1 \quad \text{for all } k \geq 1.$$

The positive solution of the quadratic equation $x^2 - 4x - 1 = 0$ is $2 + \sqrt{17}/2$. Hence the result. $\square$

**Acknowledgments.** The first version of this paper was mainly concerned with Section 5, based on the ideas exposed in the more general setting of Section 4. Thanks to two referees, the present version is much improved. The general Proposition 1, the study of mixtures of uniform distributions and the use of the submultiplicative lemma in the framework of (A2) with finite $R$ all originate from the large scope report of one of them. We are grateful for all these comments and suggestions. Other valuable ideas present in the reports have been mentioned in the paper and, in particular, in Section 8. Their full exploration would deserve further developments which are left for future work.

GET/TELECOM PARIS
CNRS LTCI
46 RUE BARRAULT
75634 PARIS CEDEX 13
FRANCE
E-MAIL: roueff@tsi.enst.fr

CENTRE FOR MATHEMATICAL SCIENCES
LUND UNIVERSITY
BOX 118
221 00 LUND
SWEDEN